\documentclass[12pt]{article}
\usepackage[utf8]{inputenc}
\usepackage{xcolor}
\usepackage[portrait, margin=1in]{geometry}
\usepackage{amsfonts}
\usepackage{amsmath}
\usepackage{amssymb}
\usepackage{amsthm}
\usepackage{enumerate}
\usepackage{comment}
\usepackage{sectsty}
\usepackage{esint}
\usepackage{graphicx}
\usepackage{float}

\newcommand{\R}{\mathbb{R}}

\newcommand{\N}{\mathbb{N}}
\newcommand{\Z}{\mathbb{Z}}

\newcommand{\B}{\mathcal{B}}

\newcommand{\CC}{\mathcal{C}}

\newcommand{\Dd}{\mathcal{D}}

\newcommand{\GG}{\mathcal{G}}

\newcommand{\Leb}{\mathcal{L}}

\newcommand{\eps}{\varepsilon}

\newcommand{\loc}{\mathrm{loc}}
\newcommand{\Om}{\Omega}

\newcommand{\bmat}{\begin{bmatrix}}
\newcommand{\emat}{\end{bmatrix}}
\newcommand{\til}{\tilde}
\newcommand{\wtil}{\widetilde}

\newcommand{\st}{\text{ s.t. }}

\DeclareMathOperator{\Lip}{Lip}

\DeclareMathOperator{\diam}{diam}
\DeclareMathOperator{\rad}{rad}

\def\bint{{\ifinner\rlap{\bf\kern.35em--}
\int\else\rlap{\bf\kern.45em--}\int\fi}\ignorespaces}

\def\bbint{{\ifinner\rlap{\bf\kern.35em--}
\hspace{0.078cm}\int\else\rlap{\bf\kern.45em--}\int\fi}\ignorespaces}

\newtheorem{theorem}{Theorem}[section]
\newtheorem{proposition}[theorem]{Proposition}

\newtheorem{corollary}[theorem]{Corollary}
\newtheorem{lemma}[theorem]{Lemma}
\newtheorem{definition}[theorem]{Definition}

\numberwithin{theorem}{section}
\numberwithin{equation}{section}

\theoremstyle{definition}
\newtheorem{remark}[theorem]{Remark}

\makeatletter
\def\blfootnote{\xdef\@thefnmark{}\@footnotetext}
\makeatother

\title{Geometric inequalities related to fractional perimeter: fractional Poincar\'e,  isoperimetric, and boxing inequalities in metric measure spaces}

\author{Josh Kline, Panu Lahti, Jiang Li, Xiaodan Zhou}
\date{\today}

\begin{document}

\blfootnote{{\it Acknowledgements:} The research of J.K. was partially supported by the University of Cincinnati Taft Research Center Dissertation Fellowship Award.  The research of J.L. was supported by the National Natural Science Foundation of China (No. 12371071 and No. 12571081). The research of X.Z. was supported by JSPS Grant-in-Aid for Scientific Research (No. 25K00211). The authors would like to thank Julian Weigt for pointing out the reference \cite{MPW24}, which inspired some of our main results, as well as anonymous referee, whose helpful comments improved the presentation of the paper.}

\maketitle

\begin{abstract}
    In the setting of a complete, doubling metric measure space $(X,d,\mu)$ supporting a $(1,1)$-Poincar\'e inequality, we show that for all $0<\theta<1$, the following fractional Poincar\'e inequality holds for all balls $B$ and locally integrable functions $u$,
    \[
    \fint_{B}|u-u_B|d\mu\le C(1-\theta)\rad(B)^\theta\fint_{\tau B}\int_{\tau B}\frac{|u(x)-u(y)|}{d(x,y)^\theta\mu(B(x,d(x,y)))}d\mu(y)d\mu(x),
    \]
    where $C\ge 1$ and $\tau\ge 1$ are constants depending only on the doubling and $(1,1)$-Poincar\'e inequality constants.
    Notably, this inequality features the scaling constant $(1-\theta)$ present in the Bourgain-Brezis-Mironescu theory characterizing Sobolev functions via nonlocal functionals.  
    
    From this inequality, we obtain a fractional relative isoperimetric inequality as well as global and local versions of a fractional boxing inequality,     
    each featuring the same scaling constant $(1-\theta)$ and defined in terms of the fractional $\theta$-perimeter, and prove equivalences with the above fractional Poincar\'e inequality.  We also show that $(X,d,\mu)$ supports a $(1,1)$-Poincar\'e inequality if and only if the above fractional Poincar\'e inequality holds for all $\theta$ sufficiently close to $1$. 

    Under the additional assumption of lower Ahlfors $Q$-regularity of the measure $\mu$, we additionally use the aforementioned results to establish global inequalities, in the form of fractional isoperimetric and fractional Sobolev inequalities, which also feature the scaling constant $(1-\theta)$.  Moreover, we prove that such inequalities are equivalent with the lower Ahlfors $Q$-regularity condition on the measure.    
\end{abstract}

\tableofcontents

\section{Introduction}

In Euclidean spaces $\mathbb{R}^n$ with $n\ge 2$, the classical isoperimetric inequality states that all bounded sets $E\subset \mathbb{R}^n$ with smooth boundary satisfy

\[
|E|^{\frac{n-1}{n}}\le \frac{1}{n|B_1|^{\frac{1}{n}}}P(E).
\]
Here, $|E|$ denotes the $n$-dimensional Lebesgue measure of $E$, $B_1$ is a unit ball in $\mathbb{R}^n$, and the perimeter of $E$, denoted as $P(E)$, coincides with the $(n-1)$-dimensional Hausdorff measure of the boundary $\partial E$ in this case.  Notably, equality holds if and only if the set $E$ is a ball. The isoperimetric inequality is shown to be equivalent to the Gagliardo-Nirenberg-Sobolev embedding inequality when $p=1$ with the sharp constant by Federer and Fleming \cite{FF60}, and Maz'ya \cite{Ma60}. The embedding inequality states that for any $f\in BV(\mathbb{R}^n)$, 
  \[
  \|f\|_{L^{\frac{n}{n-1}}(\mathbb{R}^n)} \le \frac{1}{n|B_1|^{\frac{1}{n}}} \|Df\|(\mathbb{R}^n) .
  \]

Furthermore, a local version of the isoperimetric inequality is given by the relative isoperimetric inequality. Let $B_r=B_r(x)\subset \mathbb{R}^n$ be an open ball centered at $x$ with radius $r$ and $E$ be as above, we have

\[
    \min\{|E\cap B_r|, |B_r\setminus E|\}^{\frac{n-1}{n}}\le C(n)P(E, B_r),
\]
where $P(E,B_r)$ denotes the perimeter of $E$ in $B_r$ and $C(n)$ is a constant depending only on the dimension. The relative isoperimetric inequality is known to be equivalent to the Sobolev-Poincar\'e inequality that states for any $f\in BV_{\rm loc}(\mathbb{R}^n)$, 
   \[
    \|f-f_{B_r}\|_{L^{\frac{n}{n-1}}(B_r)}\le C(n)\| Df\|(B_r).
\]
Here, $f_{B_r}=\frac{1}{|B_r|}\int_{B_r} f dx=\fint_{B_r} f dx$ and $\|Df\|$ denotes the total variation measure of $f\in BV(\mathbb{R}^n)$ \cite{EGbook, Mbook, Zbook}. 
It is not hard to see that the relative isoperimetric inequality implies a global isoperimetric inequality with a non-sharp constant when taking $B_r$ sufficiently large and the reverse also holds \cite{FFL06}. The study of the above inequalities 
plays a fundamental role in many areas, including geometry, potential theory, partial differential equations, and probability. These inequalities have been extended to Riemannian manifolds and general metric measure spaces.

 Motivated by the developments and applications in harmonic analysis, functional analysis, partial differential equations, and calculus of variations, the fractional Sobolev functions, BV functions and nonlocal fractional perimeter functionals are extensively studied in the Euclidean spaces, see  \cite{BBM01,BBM02, Br02, MaSh02, Da02, CaRoSa10, CaVa11, DiPaVa12, EEbook23, Le23} and the references therein. The fractional isoperimetric inequality and fractional Sobolev-Poincar\'e inequality and their applications have been widely explored in the literature and we list only a selection of relevant work here \cite{Po04, FS08, FuMiMo11, FiFuMaMiMo15, HuVa13, DyIhVa16, PS20, ACPS21, HuMaPeVa23, Sk23, MPW24}.

Let $\theta\in (0,1)$ and we define the fractional $\theta$-perimeter of $E$ in an open set $\Omega\subset\mathbb{R}^n$  by  \begin{equation}\label{eq:intro euc frac per}
P_\theta(E,\Omega)=\int_{E\cap\Omega}\int_{\Omega\setminus E} \frac{1}{|x-y|^{n+\theta}}dydx=\frac{1}{2}[\chi_E]_{W^{\theta, 1}(\Omega)},
\end{equation}
where $[\cdot]_{W^{\theta,1}(\Omega)}$ 
is the fractional Sobolev seminorm. We write $P_\theta(E)$ when $\Omega=\mathbb{R}^n$. The fractional isoperimetric inequality with sharp constant in the Euclidean space is proved by Frank and Seiringer \cite{FS08} that for any Borel set $E\subset \mathbb{R}^n$ with finite Lebesgue measure 
\begin{equation}\label{EGFiso}
|E|^{(n-\theta)/n}
	\le C_{n,\theta}P_\theta(E),
	\end{equation}
	 and the equality holds if and only if $E$ is a ball \cite[(4.2)]{FS08}. Recall the fractional Poincar\'e inequality in the Euclidean space shown by Bourgain, Brezis and Mironescu \cite[Theorem 1]{BBM02} states that for any cube $Q\subset \mathbb{R}^n$ and $f\in L^1_{\rm loc}(\mathbb{R}^n)$,
\begin{equation}\label{BBMFP1}
\left(\fint_Q |f-f_{Q}|^{\frac{n}{n-\theta}}dx\right)^{\frac{n-\theta}{n}}\le C(n)(1-\theta){\rm diam}(Q)^\theta\fint_Q\int_Q \frac{|f(x)-f(y)|}{|x-y|^{n+\theta}}dydx,
\end{equation}
which easily implies that
\begin{equation}\label{BBMFP}
\fint_Q |f-f_{Q}|dx\le C(n)(1-\theta){\rm diam}(Q)^\theta\fint_Q\int_Q \frac{|f(x)-f(y)|}{|x-y|^{n+\theta}}dydx.
\end{equation}

For any Borel set $E\subset \mathbb{R}^n$ with finite Lebesgue measure, one can deduce fractional relative isoperimetric inequalities by letting $f=\chi_E$ and $f=\chi_{E^c}$ that
\begin{equation}\label{eq: EFRiso1}
{\min\{ |E\cap Q|, |Q\setminus E|\}}^{\frac{n-\theta}{n}}\le C(1-\theta)P_\theta(E,Q),
\end{equation}
and
\begin{equation}\label{eq:EFRiso}
\min\{ |E\cap Q|, |Q\setminus E|\}\le C(1-\theta)\diam(Q)^\theta P_\theta(E,Q).
\end{equation}

We note that these Euclidean fractional Poincar\'e and relative isoperimetric inequalities feature the coefficient $(1-\theta)$ on the right-hand side. In the celebrated work of Bourgain, Brezis, and Mironescu \cite{BBM01}, it was shown that the Sobolev energy is explicitly recovered by multiplying this coefficient to the fractional Sobolev seminorm and letting $\theta\to 1^-$.  In particular, for the case of the fractional perimeter, one has

\begin{equation}\label{asymptotic}
\lim_{\theta\to 1^-}(1-\theta)P_\theta (E)=C(n) P(E), \quad\text{and} \quad \lim_{\theta\to 0^+}\theta P_\theta(E)=n|B_1||E|,
    \end{equation}
by the work of D\'avila \cite{Da02}, and Maz\'ya and Shaposhnikova \cite{MaSh02}.  This scaling coefficient also allows one to recover the classical Poincar\'e inequality, as in \eqref{eq:PI}, from the fractional Poincar\'e inequality by sending $\theta$ to $1$. Without the coefficient $(1-\theta)$, the right-hand side of \eqref{BBMFP1} and \eqref{BBMFP} blow up for all non-constant functions \cite{Br02}. 
 
  In the past three decades, there has been significant progress in
the study of various aspects of first-order analysis on metric measure spaces including the
theory of Sobolev spaces, BV functions, and their relation to
variational problems and partial differential equations, see \cite{Ha03, Mi03, BB11, HKST} and the references therein. In recent years, fractional Sobolev-type spaces have also undergone extensive study in the metric setting, in particular related to traces of first-order Sobolev spaces, construction of nonlocal operators, and the study of nonlocal minimal surfaces, for a sampling see \cite{GoKoSha10,Maly17,BBS22,CKKSS25,Kl25}.
In the spirit of \cite{BBM01}, characterizations of BV and Sobolev spaces via nonlocal functionals have also recently attracted great attention in the metric setting.  In particular, if $(X,d,\mu)$ is a complete, doubling metric measure space supporting a $(1,1)$-Poincar\'e inequality, see Section~\ref{sec:doubling} and Section~\ref{sec:finite per, PI}, then it was shown in \cite{DiMaS} that 
\begin{equation}\label{eq:Intro - metric BBM}
    P(E,X)\lesssim\liminf_{\theta\to 1^-}(1-\theta)P_\theta(E,X)\le\limsup_{\theta\to 1^-}(1-\theta)P_\theta(E,X)\lesssim P(E,X)
\end{equation}
for every $\mu$-measurable $E\subset X$.  Here, $P(E,X)$ is the perimeter measure of $E$ in $X$, see Section~\ref{sec:finite per, PI}, and for an open set $\Omega\subset X$ the fractional perimeter $P_\theta(E,\Omega)$ of $E$ in $\Omega$ is defined by 
\begin{equation*}
    P_\theta(E,\Omega):=\int_{\Omega\cap E}\int_{\Omega\setminus E}\frac{2}{d(x,y)^\theta[\mu(B(y,d(x,y)))+\mu(B(x,d(x,y)))]}d\mu(y)d\mu(x).
\end{equation*}
For generalizations of this result, see also \cite{Go22, Ha24, LPZ, LPZ2, Shpre}. From \eqref{eq:Intro - metric BBM}, we see that in this setting, the $\theta$-perimeter recovers the classical perimeter, up to comparability, as $\theta\to 1^-$, provided one rescales by the crucial $(1-\theta)$ coefficient.  Since in the Euclidean setting, the fractional Poincar\'e inequalities \eqref{BBMFP} and \eqref{BBMFP1}, and consequently the fractional relative isoperimetric inequalities \eqref{eq:EFRiso} and \eqref{eq: EFRiso1}, hold with such a scaling coefficient \cite{BBM02}, it is natural to ask whether such Bourgain-Brezis-Mironescu-type fractional inequalities also hold in the metric setting.  The main goal of this paper is to establish these inequalities and explore the relationships between them in the setting of a complete, doubling metric measure space supporting a $(1,1)$-Poincar\'e inequality.  The notions of a doubling measure and $(1,1)$-Poincar\'e inequality will be introduced in Sections~\ref{sec:doubling} and \ref{sec:finite per, PI}, respectively. This setting encompasses a wide variety of spaces including Euclidean and Muckenhoupt weighted Euclidean spaces, a large class of domains in $\mathbb{R}^n$ satisfying appropriate conditions, Riemannian manifolds with non-negative Ricci curvature, Carnot groups, and Laakso spaces, for example.   

The first of our main results is the following, which generalizes \eqref{BBMFP1} and \eqref{BBMFP}:

\begin{theorem}\label{thm: MFPinq}
  Let  $(X,d,\mu)$ be a complete, doubling metric measure space supporting a $(1,1)$-Poincar\'e inequality.  Let $0<\theta<1$ and $1\le q\le Q_d/(Q_d-\theta)$, where $Q_d>1$ is the relative lower mass bound exponent \eqref{eq:rel lower mass bound exponent}.  Then there exist constants $C\ge 1$ and $\tau\ge 1$ so that for all balls $B:=B(x,r)\subset X$ and all $u\in L^1_{\rm loc}(X)$, we have 
    \begin{equation}\label{MFPinq}
    \left(\fint_{B}|u-u_{B}|^qd\mu\right)^{1/q}\le C(1-\theta)r^\theta\fint_{\tau B}\int_{\tau B}\frac{|u(x)-u(y)|}{d(x,y)^\theta\mu(B(x,d(x,y)))}d\mu(y)d\mu(x).
    \end{equation}
    Here, the constants $C$ and $\tau$ depend only on the doubling and $(1,1)$-Poincar\'e inequality constants.
\end{theorem}

If \eqref{MFPinq} holds, then we say that $(X,d,\mu)$ supports a \emph{ $(\theta,q,1)_{BBM}$-Poincar\'e inequality}, see Definition~\ref{def:frac ineq}.  We use this notation to emphasize the presence of the $(1-\theta)$ coefficient and to differentiate this inequality from versions without, such as those in \cite{DyLeVa23,CeDrMa19}. Our proof of the above theorem is inspired by the recent work by Myyryl\"ainen, P\'erez, and Weigt \cite{MPW24}.  In particular, we establish the above result by first proving an inequality in Lemma \ref{lem:Fractional Isoperimetric} which acts as an improvement of the standard relative isoperimetric inequality.  See Remark \ref{improvedRiso} for more details.

As a complete, doubling metric measure space supporting a $(1,1)$-Poincar\'e inequality admits a geodesic metric which is bilipschitz equivalent with the original metric, see for example \cite[Corollary 8.3.16]{HKST}, it suffices to prove Theorem~\ref{thm: MFPinq} under the assumption that $(X,d)$ is geodesic, see the discussion at the beginning of Section~\ref{FPI}.  As the proof of Theorem~\ref{thm: MFPinq} shows, we may take $\tau=1$ in \eqref{MFPinq} when $(X,d,\mu)$ is additionally assumed to be geodesic.  This mirrors the classical result for complete, doubling metric measure spaces supporting a $(1,p)$-Poincar\'e inequality \cite{HaKo00}.

As a direct application of the fractional Poincar\'e inequality \eqref{MFPinq}, we obtain a fractional relative isoperimetric inequality: 

\begin{corollary}\label{thm: MFRiso}
   Let  $(X,d,\mu)$ be a complete, doubling metric measure space supporting a $(1,1)$-Poincar\'e inequality.  Let $0<\theta<1$ and $1\le q\le Q_d/(Q_d-\theta)$, where $Q_d>1$ is the relative lower mass bound exponent \eqref{eq:rel lower mass bound exponent}.  Then there exist constants $C\ge 1$ and $\tau\ge 1$, so that for all measurable sets $E\subset X$, and all balls $B:=B(x,r)\subset X$, we have that 
    \begin{equation}\label{MFRiso}
        \left(\frac{\min\{\mu(B\cap E),\,\mu(B\setminus E)\}}{\mu(B)}\right)^{1/q}\le C(1-\theta)r^\theta\frac{ P_\theta(E,\tau B)}{\mu(\tau B)}.
    \end{equation}
    Here the constants $C$ and $\tau$ depend only on the doubling and $(1,1)$-Poincar\'e inequality constants.
\end{corollary}

 If \eqref{MFRiso} holds, then we say that $(X,d,\mu)$ supports a \emph{$(\theta,q,1)_{BBM}$-relative isoperimetric inequality}, see Definition~\ref{def:frac ineq}. There are various results of isoperimetric inequalities \cite{We05, KoLa14, CaMo17, Na22} in metric spaces, but to our best knowledge, the fractional isoperimetric inequality has not been previously investigated in the general setting.

In \cite{KoLa14}, it was shown that a doubling metric measure space supports a $(q,1)$-Poincar\'e inequality if and only if it supports a $(q,1)$-relative isoperimetric inequality, see \eqref{eq:rel isoperimetric}.  We show here that the analogous result holds in the fractional case. 

\begin{theorem}\label{rmk: PiRiequiv}
    Let $0<\theta<1$ and $q\ge 1$.  Then $(X,d,\mu)$ supports a $(\theta,q,1)_{BBM}$-Poincar\'e inequality if and only if it supports a $(\theta,q,1)_{BBM}$-relative isoperimetric inequality. 
\end{theorem}

The assumption of a (non-fractional) $(1,p)$-Poincar\'e inequality imposes significant geometric restrictions on a metric measure space $(X,d,\mu)$.  For example, such spaces are necessarily quasiconvex, see \cite{HKST}, and for any $p_0\ge 1$, there are examples of metric measure spaces which support $(1,p)$-Poincar\'e inequalities for all $p>p_0$, but not $p=p_0$. See \cite{MTW13} or for a simple example, consider the bowtie spaces.  In contrast, fractional Poincar\'e inequalities, without the $(1-\theta)$ coefficient, do not carry such geometric information and hold under rather mild assumptions on the space.  For example, an analog of the fractional Poincar\'e inequality \eqref{BBMFP}, without the $(1-\theta)$ coefficient, follows immediately if $\mu$ is doubling, see \cite[Lemma~2.2]{DyLeVa23}. Furthermore, using a Maz'ya-type truncation argument, it was also shown in \cite[Theorem~3.4]{DyLeVa23} that an analog of \eqref{BBMFP1}, with $n$ replaced by the relative lower mass bound exponent from \eqref{eq:rel lower mass bound exponent} and without the constant $(1-\theta)$, holds when $\mu$ is both doubling and reverse doubling, see Section~\ref{sec:doubling}.  See also \cite{CeDrMa19} for formulations of this inequality on John domains in the metric measure space setting.  In order to obtain the improved versions of these inequalities in Theorem~\ref{thm: MFPinq}, we additionally assume that $(X,d,\mu)$ supports a $(1,1)$-Poincar\'e inequality.  In fact, this assumption is necessary to obtain the $(1-\theta)$ coefficients in the fractional Poincar\'e inequality:

\begin{theorem}\label{fractional to Poincare}
	Let  $(X,d,\mu)$ be a complete, doubling metric measure space and let $1\le q<\infty$.  Suppose
	there exist constants $C,\tau \ge 1$ and a sequence of numbers $\theta_i\in (0,1)$ converging to $1$, such that for all balls $B:=B(x,r)\subset X$, and all $u\in L^1_\loc(X)$, we have 
	\[
	\left(\fint_{B}|u-u_{B}|^q\,d\mu\right)^{1/q}
	\le C(1-\theta_i)r^{\theta_i}\fint_{\tau B}\int_{\tau B}\frac{|u(x)-u(y)|}{d(x,y)^{\theta_i}\mu(B(x,d(x,y)))}\,d\mu(x)\,d\mu(y).
	\]
	Then $(X,d,\mu)$ supports a $(q,1)$-Poincar\'e inequality.
\end{theorem}

In this manner, the $(\theta,q,1)_{BBM}$-Poincar\'e inequalities encode the same geometric information as the $(1,1)$-Poincar\'e inequality. Moreover, since a $(q,1)$-Poincar\'e inequality is equivalent to a $(q,1)$-relative isoperimetric inequality by \cite{KoLa14}, and the same holds for their fractional counterparts by Theorem~\ref{rmk: PiRiequiv} above, it follows by Theorem~\ref{fractional to Poincare} that the $(q,1)$-relative isoperimetric inequality is recovered by its fractional counterpart.  That is, if a complete doubling metric measure space supports a $(\theta,q,1)$-relative isoperimetric inequality for all $1-\eps<\theta<1$, then it supports a $(q,1)$-relative isoperimetric inequality. Theorem~\ref{fractional to Poincare} is established in Section~\ref{FPIeqv} as a corollary to Theorem~\ref{thm:rho and Poincare}, which is proven for more general kernels in the fractional Poincar\'e inequality.

Next, we apply the fractional Poincar\'e inequality \eqref{MFPinq} to deduce a fractional boxing inequality. The (non-fractional) boxing inequality in the Euclidean setting, first proved by Gustin \cite{Gu60}, states that every open bounded set $U\subset \mathbb{R}^n$ can be covered by a collection of balls $U\subset \bigcup_{i=1}^\infty B_{r_i}$ such that $\sum_{i=1}^\infty r_i^{n-1}\le C(n)P(U)$, where $C(n)$ depends only on the dimension. Such a boxing inequality was generalized to the setting of doubling metric measure spaces supporting a $(1,1)$-Poincar\'e inequality by Kinnunen, Korte, Shanmugalingam and Tuominen \cite{KKST2}.  For more on the extension of the boxing inequality to Banach spaces and Riemannian manifolds, see the recent work \cite{AvNa25}.

In the Euclidean setting, Ponce and Spector \cite{PS20} established the following fractional boxing inequality, involving the $\theta$-perimeter: given an open bounded set $U\subset\R^n$, one can find a collection of balls covering $U$ such that
\[
\sum_{i=1}^\infty r_i^{n-\theta}\le C(n)\theta(1-\theta)P_\theta (U),
\]
for every $\theta\in (0,1)$. As $\theta\to 1^-$, one recovers the Gustin boxing inequality.  Applying the $(\theta,1,1)_{BBM}$-Poincar\'e inequality obtained in Theorem~\ref{thm: MFPinq}, we deduce the following fractional boxing inequality in the metric setting, using ideas from the Euclidean proof in \cite{PS20} and the proof in \cite{KKST2} of the non-fractional boxing inequality in the metric setting: 

\begin{theorem}\label{MFBoxing}
    Let  $(X,d,\mu)$ be a complete, doubling metric measure space supporting a $(1,1)$-Poincar\'e inequality, such that $\mu(X)=\infty$.  Then there exist constants $C\ge 1$ and $\tau\ge 1$ such that for any bounded measurable set  $U\subset X$, there exists a set $\mathcal{N}\subset U$ with $\mu(\mathcal{N})=0$ and a countable collection of disjoint balls $\{B(x_i,r_i)\}_{i\in\mathbb{N}}$ such that 
    \[
U\setminus \mathcal{N}\subset\bigcup_{i\in\mathbb{N}}B(x_i,5\tau r_i), \quad\quad\frac{1}{2C_\mu}\le\frac{\mu(B(x_i,r_i)\cap U)}{\mu(B(x_i,r_i))}<\frac{1}{2}
    \]
    for each $i\in\mathbb{N}$, and for each $0<\theta<1$, we have 
\begin{equation}\label{eq:doubling boxing ineq}
        \sum_{i\in\mathbb{N}}\frac{\mu(B(x_i,r_i))}{r_i^\theta}\le C\theta(1-\theta)P_\theta(U,X).
    \end{equation}
    Here, the constants $C$ and $\tau$ depend only on the doubling and $(1,1)$-Poincar\'e inequality constants.
\end{theorem}
If $U$ is additionally assumed to be open, then the conclusion of the above theorem holds without requiring the negligible set $\mathcal{N}$, see Remark~\ref{rem:open set boxing}.

By using the $(\theta,1,1)_{BBM}$-relative isoperimetric inequality, obtained in Corollary~\ref{thm: MFRiso}, we also establish in Proposition~\ref{thm:loc-box-ineq} a local version of the fractional boxing inequality above.  In this result, $U$ may be taken to be an arbitrary measurable subset of a given ball $B$, contained up to a set of measure zero in the union of the  desired collection of inflated balls, and the fractional perimeter in \eqref{eq:doubling boxing ineq} is given by $P_\theta(U,B)$.  This local fractional boxing inequality in turn implies the $(\theta,q,1)_{BBM}$-Poincar\'e inequality, as shown in Corollary~\ref{cor:local box implies frac PI}.

In order to obtain a fractional isoperimetric inequality generalizing \eqref{EGFiso}, as well as an equivalent fractional Sobolev-type embedding, we need a further assumption on the measure $\mu$. Namely, we assume the measure $\mu$ is \emph{lower Ahlfors $Q$-regular}, that is, there exist constants $c_0>0$ and $Q>0$ such that for every $x\in X$ and $0<r<2\diam (X)$, we have
\begin{equation}\label{eq:lower Q}
\mu(B(x,r))\ge c_0 r^Q.
\end{equation}
Note that when $X$ is bounded and $\mu$ is doubling, the measure $\mu$ is lower Ahlfors $Q_d$-regular, where $Q_d$ is the relative lower mass bound exponent of $\mu$, see \eqref{eq:rel lower mass bound exponent}.

By applying Theorem~\ref{MFBoxing} with the additional assumption of lower Ahlfors $Q$-regularity, we obtain the following two inequalities, which we refer to as a \emph{$(\theta, \frac{Q}{Q-\theta},1)_{BBM}$-isoperimetric inequality} and a \emph{$(\theta, \frac{Q}{Q-\theta},1)_{BBM}$-Sobolev inequality} respectively, see Definition~\ref{def:frac ineq}:

\begin{corollary}\label{cor:Global inequalities}
    Let $(X,d,\mu)$ be a complete, doubling metric measure space supporting a $(1,1)$-Poincar\'e inequality.  Assume that $\mu(X)=\infty$ and that $\mu$ satisfies the lower Ahlfors $Q$-regularity condition \eqref{eq:lower Q} for $Q>1$.  Then, 
    \begin{enumerate}
    \item [(i)]  there exists $C\ge 1$, depending only on the doubling, lower Ahlfors $Q$-regularity, and $(1,1)$-Poincar\'e inequality constants such that for all $0<\theta<1$ and all $\mu$-measurable $E\subset X$ such that $\mu(E)<\infty$, we have \begin{equation}\label{MFGiso}
	\mu(E)^{(Q-\theta)/Q}
	\le C\theta(1-\theta)P_\theta(E,X).
	\end{equation}
    \item[(ii)]  there exists $C\ge 1$, depending only on the doubling, lower Ahlfors $Q$-regularity, and $(1,1)$-Poincar\'e inequality constants such that for all $0<\theta<1$ and all $u\in L^1(X)$, we have
     \begin{equation}\label{GFPinq}
        \left(\int_{X}|u|^{\frac{Q}{Q-\theta}}d\mu\right) ^{\frac{Q-\theta}{Q}}\le C\theta(1-\theta)\int_{X}\int_{X}\frac{|u(x)-u(y)|}{d(x,y)^\theta\mu(B(x,d(x,y)))}d\mu(y)d\mu(x).
    \end{equation}
    \end{enumerate}
    Furthermore $(X,d,\mu)$ supports a $(\theta, \frac{Q}{Q-\theta},1)_{BBM}$-isoperimetric inequality if and only if it supports a $(\theta, \frac{Q}{Q-\theta},1)_{BBM}$-Sobolev inequality.
\end{corollary}

\begin{remark}
Based on the asymptotic relation of fractional perimeter \eqref{asymptotic}, it is clear that the order of $\theta$ is sharp in our Corollary~\ref{cor:Global inequalities}. Note the sharp constant of this inequality proved by Frank and Seiringer in the Euclidean space \cite[(4.2)]{FS08} is given by an integral formula and the order of $\theta$ near $0$ and $1$ is less obvious compared to our theorem. Then letting $\theta \rightarrow 1$ in \eqref{MFGiso}, by \eqref{eq:Intro - metric BBM}, we have  isoperimetric inequality in a metric measure space, that is, for every $\mu$-measurable $E\subset X$ with $\mu(E)<\infty$, one has
\begin{eqnarray}\label{eq:iso-ineq}
    \mu(E)^{(Q-1)/Q}
	\le CP(E,X),
\end{eqnarray}
where $C$ only depends on the doubling, lower Ahlfors $Q$-regularity, and $(1,1)$-Poincar\'e inequality constants. Moreover, \eqref{eq:iso-ineq} is the main result of Jiang and Koskela \cite{JK12}.  Under the conditions of Corollary \ref{cor:Global inequalities}, we find that $(X,d,\mu)$ satisfying
local quantitative Lipschitz regularity for solutions to the Cheeger-Poisson equation, introduced by Jiang and Koskela, is not essential. For more details see \cite{JK12}.
\end{remark}

If $\mu$ is additionally assumed to be lower Ahlfors $Q$-regular, then as an application of the local fractional boxing inequality Proposition~\ref{thm:loc-box-ineq}, we also obtain versions of the $(\theta,q,1)_{BBM}$-Poincar\'e and relative isoperimetric inequalities, with $q=Q/(Q-\theta)$, without the measure of the ball and radius terms, see Corollary~\ref{cor:improved inequalities}.  We refer to these inequalities respectively as \emph{improved $(\theta,\frac{Q}{Q-\theta},1)_{BBM}$-Poincar\'e} and \emph{improved $(\theta,\frac{Q}{Q-\theta},1)_{BBM}$-relative isoperimetric inequalities}, see Definition~\ref{def:frac ineq}. Furthermore, these improved inequalities are also equivalent, by the arguments of Theorem~\ref{rmk: PiRiequiv}.

It turns out that the lower Ahlfors $Q$-regularity condition is essential to obtain the $(\theta,\frac{Q}{Q-\theta},1)_{BBM}$-Sobolev inequality or the improved $(\theta,\frac{Q}{Q-\theta},1)_{BBM}$-Poincar\'e inequality (or their equivalent isoperimetric inequalities), as shown in the following theorem. This result is analogous to the results showing the equivalence of lower Ahlfors $Q$-regularity and the validity of Sobolev embedding theorems in Euclidean domains \cite{HKT08, Z15}, in metric spaces \cite{AGH20} and in quasi-metric spaces \cite{AYY24}. 
\begin{theorem}\label{lowerQEquiv}
    Let  $(X,d,\mu)$ be a complete, doubling metric measure space supporting a $(1,1)$-Poincar\'e inequality such that $\mu(X)=\infty$.  Then, the following conditions are equivalent.
    \begin{enumerate}
        \item[(i)] $\mu$ satisfies the lower Ahlfors $Q$-regularity condition \eqref{eq:lower Q} for $Q>1$. 
        \item[(ii)] There exists a constant $C\ge 1$ such that for all $0<\theta<1$ and $u\in L^1(X)$, we have
\[
\left(\int_{X}|u|^{\frac{Q}{Q-\theta}}d\mu\right) ^{\frac{Q-\theta}{Q}}\le C\theta(1-\theta)\int_{X}\int_{X}\frac{|u(x)-u(y)|}{d(x,y)^\theta\mu(B(x,d(x,y)))}d\mu(y)d\mu(x).
\]
    \item[(iii)] There exist constants $C\ge 1$ and $\tau\ge 1$ such that for all $0<\theta<1$, and all balls $B_0:=B(x_0,r_0)\subset X$ and $u\in L^1(X)$, we have
    \[
    \left(\int_{B_0}|u-u_{B_0}|^{\frac{Q}{Q-\theta}}d\mu\right) ^{\frac{Q-\theta}{Q}}\le C(1-\theta)\int_{\tau B_0}\int_{\tau B_0}\frac{|u(x)-u(y)|}{d(x,y)^\theta\mu(B(x,d(x,y)))}d\mu(y)d\mu(x).
    \]
    \end{enumerate}
The lower Ahlfors $Q$-regularity constants in the implications $(ii)\implies (i)$ and $(iii)\implies (i)$ are independent of $\theta$, while  the constants $C$ and $\tau$ in $(ii)$ and $(iii)$ may depend additionally on the doubling and $(1,1)$-Poincar\'e inequality constants for the implications $(i)\implies (ii)$ and $(i)\implies (iii)$. 
\end{theorem}

A diagram of the main results proved in our paper on a complete doubling metric measure space is given below. The orange boxes contain the underlying conditions on the spaces, the green boxes contain functional inequalities while the red boxes include geometric inequalities. 
\begin{figure}[H]
\centering
  \includegraphics[width=7\textwidth,height=11cm,keepaspectratio]{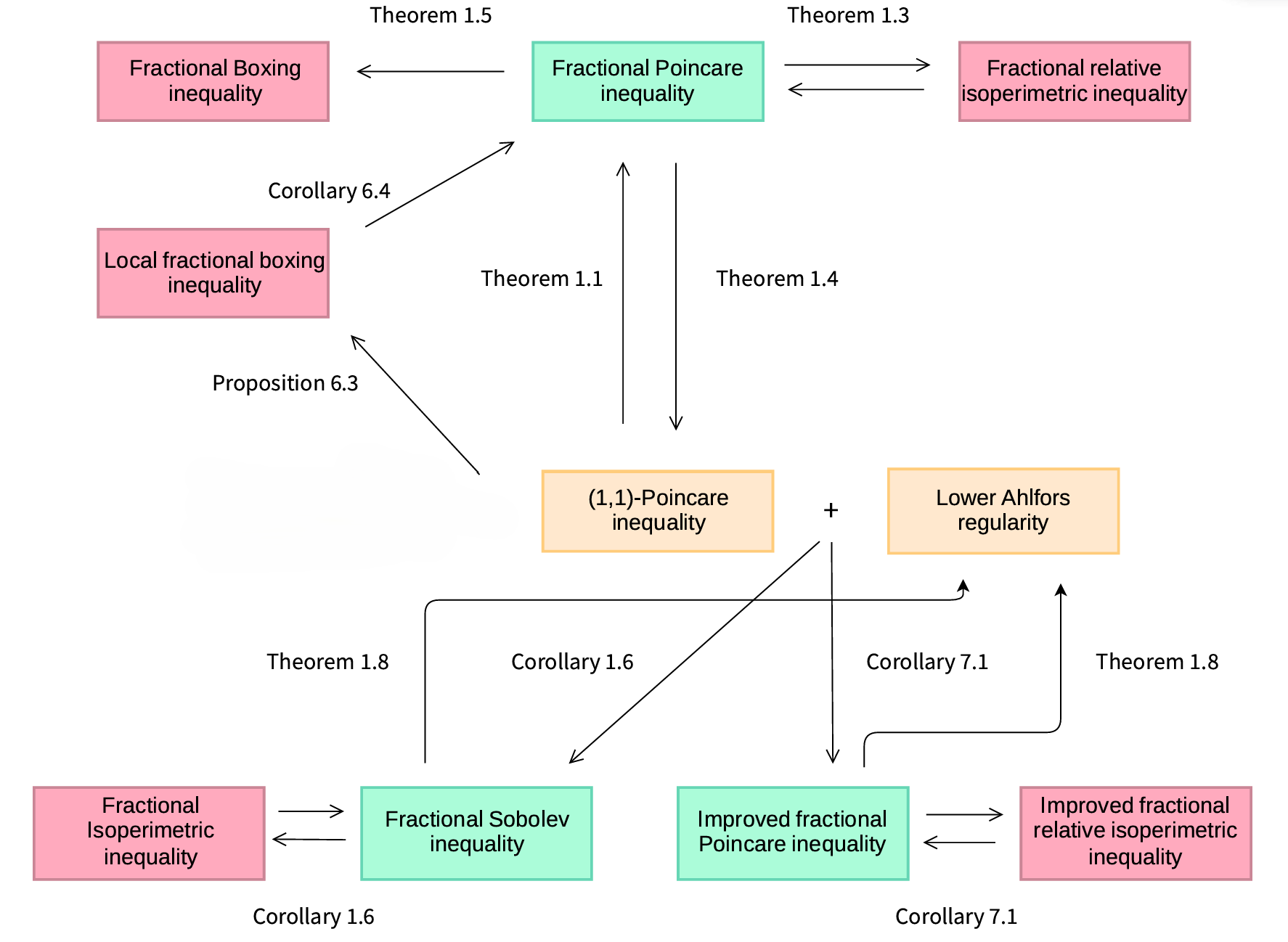}
\caption{Inequality implications}
\end{figure}

The paper is organized as follows. We will review some background notions, definitions, and show some preliminary results in Section \ref{Preliminaries}. In Section \ref{FRIIne}, we show several key lemmas which can be seen as an improved, nonlocal version of the classical relative isoperimetric inequality \eqref{eq:rel isoperimetric}. In Section \ref{FPI}, we apply the key lemmas to show the fractional Poincar\'e inequality Theorem \ref{thm: MFPinq}, fractional relative isoperimetric inequality Corollary~\ref{thm: MFRiso}, and the equivalence Theorem~\ref{rmk: PiRiequiv}. The relationship between the fractional Poincar\'e inequality and $(1,1)$-Poincar\'e inequality given by Theorem \ref{fractional to Poincare} will be established in Section \ref{FPIeqv}. The fractional boxing inequality Theorem \ref{MFBoxing} is proved in Section \ref{FBoxing}. Finally, we show the fractional isoperimetric and fractional Sobolev inequality Corollary~\ref{cor:Global inequalities}, and the equivalence with the lower Ahlfors $Q$-regular condition Theorem \ref{lowerQEquiv} in Section \ref{Fiso}.

\section{Preliminaries}\label{Preliminaries}

\subsection{Doubling measures and length spaces}\label{sec:doubling}
Throughout this paper, we assume that $(X,d,\mu)$ is a complete metric space equipped with a Borel regular outer measure $\mu$.  We also assume that $\mu$ is \emph{doubling}, that is, there exists a constant $C_\mu\ge 1$ such that for every $x\in X$ and $r>0$, we have
\[
0<\mu(B(x,2r))\le C_\mu\,\mu(B(x,r))<\infty.
\]
By iterating the doubling condition, there exist constants $C\ge 1$ and $Q_d>1$, depending only on $C_\mu$, such that 
\begin{equation}\label{eq:rel lower mass bound exponent}
\frac{\mu(B(y,r))}{\mu(B(x,R))}\ge C^{-1}\left(\frac{r}{R}\right)^{Q_d}
\end{equation}
for all $x\in X$, $0<r\le R<\infty$, and $y\in B(x,R)$ \cite[Lemma 4.7]{Ha03}.  The constant $Q_d$ is referred to as the \emph{relative lower mass bound exponent} of $\mu$. 

A measure $\mu$ is said to be \emph{reverse doubling} if there exists $s>0$ and $C_s\ge 1$,  
\begin{equation}\label{eq:intro reverse doubling}
\frac{\mu(B(y,r))}{\mu(B(x,R))}\le C_s\left(\frac{r}{R}\right)^s
\end{equation}
for all $0<r\le R\le 2\diam X$, $x\in X$ and $y\in B(x,R)$.  If $(X,d,\mu)$ is connected with $\mu$ doubling, then $\mu$ is also reverse doubling with $s$ and $C_s$ depending only on $C_\mu$, see for example \cite[Corollary~3.8]{BB11}.

Throughout this paper, we at times let $C>0$ denote a constant which, unless otherwise specified, depends only on the structural constants of the metric measure space, such as the doubling constant, or Poincar\'e inequality constants (see below), for example.  Its precise value is not of interest to us, and may change with each occurrence, even within the same line.  Furthermore, given quantities $A$ and $B$, we will often use the notation $A\simeq B$ to mean that there exists a constant $C\ge 1$ such that $C^{-1} A\le B\le CA$. Likewise, we use $A\lesssim B$ and $A\gtrsim B$ if the left and right inequalities hold, respectively.

A metric space $(X,d)$ is said to be \emph{geodesic} if for any two points $x,y\in X$, there exists a curve $\gamma:[a,b]\to X$ joining $x$ to $y$ such that $\ell(\gamma)=d(x,y)$, where $\ell(\gamma)$ is the length of $\gamma$ with respect to $d$.  A metric space $(X,d)$ is said to be a \emph{length space} if for all $x,y\in X$, we have that
\[
d(x,y)=\inf_\gamma\ell(\gamma),
\]
where the infimum is over all curves joining $x$ to $y$. 

Geodesic metric spaces equipped with a doubling measure satisfy the following Calderon-Zygmund-type decomposition for sets.  This lemma is a generalization to the metric setting of \cite[Lemma~2.1]{MPW24}.

\begin{lemma}\label{lem:the decomposition of the ball}
	Let $(X,d,\mu)$ be a geodesic metric space, with $\mu$ a doubling measure. Let $B_0:=B(x_0,r_0)\subset X$ and $E\subset X$ be a $\mu$-measurable set. Assume that
	\[
    \frac{\mu(B_0 \cap E )}{\mu(B_0)}\leq \lambda 
	\]
	holds for some $0<\lambda<1.$ Then there exist countably many pairwise disjoint balls $B_i \subset B_0$ , $i \in \mathbb{N}$ such that 
	\begin{align*}
    &(i)\quad\rad(B_i)=2^{-N_i}r_0\text{ for some }N_i\in\N\cup\{0\},\\
	&(ii)\quad B_0 \cap E \subset \bigcup_{i}5B_i \quad \textrm{up to a set of }\mu\textrm{-measure zero},\\	
	&(iii)\quad \lambda/C_\mu^2 \le\frac{\mu(B_i\cap E)}{\mu(B_i)}\leq\lambda,
	\end{align*}
	where $C_\mu$ is the doubling constant of $\mu$.
\end{lemma}

\begin{proof}
Let $B_0:=B(x_0,r_0)\subset X$.  For any ball $B:=B(x,r)\subset X$ and $k\in\N$, we define the collection of balls
\[
\Dd_k(B):=\{B'\subset B:\rad(B')=2^{-k}r_0\}.
\]
We claim that for all $k\in\N$ 
\begin{equation}\label{eq:covering B_0}
    B_0=\bigcup_{B\in\Dd_k(B_0)}B.
\end{equation}
To prove this claim, it suffices by induction to prove the case $k=1$.  To this end, let $x\in B_0$ and without loss of generality, assume that $d(x,x_0)\ge r_0/2$. Let $\delta:=r_0-d(x,x_0)\le r_0/2$.  Connecting $x$ to $x_0$ by a geodesic, there exists $y$ along this geodesic such that $d(x_0,y)=(r_0-\delta)/2$. By the geodesic property, it follows that 
\[
d(x,y)=d(x,x_0)-d(x_0,y)=d(x,x_0)-(r_0-\delta)/2=d(x,x_0)/2<r_0/2.
\]
Thus, $x\in B(y,r_0/2)$.  Moreover, if $z\in B(y,r_0/2)$, then by the triangle inequality, it follows that 
\[
d(z,x_0)\le d(z,y)+d(y,x_0)<r_0/2+(r_0-\delta)/2=r_0-\delta/2<r_0,
\]
 and so $B(y,r_0/2)\in\Dd_1(B_0)$.  This proves the claim.

We now assume that $\mu(B_0\cap E)<C_\mu^{-2}\lambda\mu(B_0)$, otherwise the collection consisting of only the ball $B_0$ satisfies the conclusion of the lemma. We define the following collection of balls
\[
\GG_1(B_0):=\{B\in\Dd_1(B_0):\mu(B\cap E)\ge C_\mu^{-2}\lambda \mu(B)\},\qquad\B_1(B_0):=\Dd_1(B_0)\setminus\GG_1(B_0),
\]
and set $\GG_1:=\GG_1(B_0)$ and $\B_1:=\B_1(B_0)$.  As $\mu$ is doubling, it follows that for each $B\in\GG_1$,
\[
\mu(B\cap E)\le\mu(B_0\cap E)\le C_\mu^{-2}\lambda\mu(B_0)\le C_\mu^{-2}\lambda\mu(4B)\le\lambda\mu(B).
\]

Inductively, having defined $\GG_k$ and $\B_k$, we define the following collection of balls  for each $B\in \B_k$:
\[
\GG_{k+1}(B):=\{B'\in\Dd_{k+1}(B):\mu(B'\cap E)\ge C_\mu^{-2}\lambda \mu(B')\},\qquad \B_{k+1}(B):=\Dd_{k+1}(B)\setminus\GG_{k+1}(B).
\]
We then set $\GG_{k+1}:=\bigcup_{B\in\B_{k}}\GG_{k+1}(B)$ and $\B_{k+1}:=\bigcup_{B\in\B_{k}}\B_{k+1}(B)$. If $B'\in\GG_{k+1}$, then there exists $B\in\B_k$ such that $B'\subset B$.  By doubling, it then follows that 
\begin{align}\label{eq:density upper bound}
\mu(B'\cap E)\le\mu(B\cap E)\le C_\mu^{-2}\lambda\mu(B)\le\lambda\mu(B').
\end{align}

Let $\GG:=\bigcup_{k\in\N}\GG_k$.  We claim that $B_0\cap E\subset\bigcup_{B\in\GG}B$ up to a set of $\mu$-measure zero.  Indeed, if $x\in (B_0\cap E)\setminus\bigcup_{B\in\GG}B$, then inductively using the claim \eqref{eq:covering B_0} above, it follows that for each $k\in\N$, there exists $B_k\in\B_k$ such that $x\in B_k$. Since  $\mu(B_k\cap E)<C_\mu^{-2}\lambda\mu(B_k)$, it follows that
\[
\left|\fint_{B_k}\chi_Ed\mu-1\right|=1-\frac{\mu(E\cap B_k)}{\mu(B_k)}>1-\frac{1}{C_\mu^2}
\]
for every $k$. On the other hand, by the fact that $\mu$ is a doubling measure, one has
\[
\left|\fint_{B_k}\chi_Ed\mu-1\right|\lesssim\fint_{B(x,2^{-k+1}r_0)}|\chi_E-1|d\mu,
\]
and the right hand side goes to zero as $k\to \infty$ whenever $x$ is a Lebesgue point of $E$ by Lebesgue differentiation theorem. This implies that $(B_0\cap E)\setminus\bigcup_{B\in\GG}B$ belongs to a null set and the claim is proved.

 Therefore, by the $5$-covering lemma, there exists a pairwise disjoint, countable subcollection $\{B_i\}_{i\in\N}\subset\GG$ such that 
\[
B_0\cap E\subset\bigcup_{i\in\N}5B_i,
\]
up to a set of $\mu$-measure zero, which gives us $(ii)$.  Finally, the collection $\{B_i\}_{i\in\N}$ satisfies $(i)$ and $(iii)$ by construction and \eqref{eq:density upper bound}.
\end{proof}

 We note that if $E$ is an open set, then (ii) holds in the above lemma without excluding a $\mu$-negligible set.

Length spaces equipped with a doubling measure satisfy the following annular decay property:

\begin{lemma}\label{lem:annular decay}\cite[Proposition~11.5.3]{HKST} Let $(X,d)$ be a length space equipped with a doubling measure $\mu$.  Then there exist constants $C_A\ge 1$ and $0<\beta\le 1$, depending only on $C_\mu$, such that 
\[
\mu(B(x,r)\setminus\overline B(x,(1-\eps)r))\le C_A\eps^\beta\mu(B(x,r))
\]
for every $x\in X$, $r>0$, and $0<\eps\le 1$.    
\end{lemma}

\subsection{Sets of finite perimeter, Poincar\'e inequalities, and isoperimetric inequalities}\label{sec:finite per, PI}
A non-negative Borel function $g:X\to[0,\infty]$ is an \emph{upper gradient} of a function $u:X\to\overline\R$ if the following holds for all non-constant, compact, rectifiable curves $\gamma:[a,b]\to X$:
\[
|u(\gamma(b))-u(\gamma(a))|\le\int_\gamma g\,ds,
\]
whenever $u(\gamma(a))$ and $u(\gamma(b))$ are both finite, and $\int_\gamma g\,ds=\infty$ otherwise. 

Complete metric spaces equipped with a doubling measure are proper, that is, closed and bounded subsets are compact.  Thus, for any open set $\Omega\subset X$, we define $\Lip_\loc(\Omega)$ to be the space of functions which are Lipschitz in any open $\Omega'\Subset\Omega$.  By $\Omega'\Subset\Omega$, we mean that $\overline\Omega'$ is a compact subset of $\Omega$.  We define other local spaces of functions similarly. 

For an open set $\Omega\subset X$ and $u\in L^1_\loc(\Omega)$, we define the \emph{total variation} of $u$ in $\Omega$ by
\[
\|Du\|(\Omega):=\inf\left\{\liminf_{i\to\infty}\int_\Omega g_{u_i}d\mu:u_i\in\Lip_\loc(\Omega),\, u_i\to u\text{ in }L^1_\loc(\Omega)\right\},
\]
where each $g_{u_i}$ is an upper gradient of $u_i$. For an arbitrary set $A\subset X$, we define 
\begin{equation}\label{eq: BV energy general sets}
\|Du\|(A):=\inf\{\|Du\|(\Omega):\Omega\subset X\text{ open},\,A\subset\Omega\}.
\end{equation}
A $\mu$-measurable set $E\subset X$ is said to be a \emph{set of finite perimeter} if $\|D\chi_E\|(X)<\infty$.  The perimeter of $E$ in $\Omega$ is denoted by 
\[
P(E,\Omega):=\|D\chi_E\|(\Omega).
\]
If $P(E,\Omega)<\infty$, then $P(E,\cdot)$ is a finite Radon measure on $\Omega$ by \cite[Theorem~3.4]{Mi03}.

For $1\le p,q<\infty$, we say that $(X,d,\mu)$ supports a \emph{$(q,p)$-Poincar\'e inequality} if there exist constants $C,\tau\ge 1$ such that for every ball $B\subset X$, every $u\in L^1_\loc(X)$, and every upper gradient $g$ of $u$, we have
\begin{equation}\label{eq:PI}
    \left(\fint_{B}|u-u_B|^q\,d\mu\right)^{1/q}\le C\rad(B)\left(\fint_{\tau B}g^p\,d\mu\right)^{1/p}.
\end{equation}

For $1\le q<\infty$, we say that $(X,d,\mu)$ supports a \emph{$(q,1)$-relative isoperimetric inequality} if there exist constants $C,\tau\ge 1$ such that for all balls $B\subset X$ and $\mu$-measurable sets $E\subset X$, we have 
\begin{equation}\label{eq:rel isoperimetric}
    \left(\frac{\min\{\mu(B\cap E),\,\mu(B\setminus E)\}}{\mu(B)}\right)^{1/q}\le C\rad(B)\frac{ P
    (E,\tau B)}{\mu(\tau B)}.
\end{equation}
By \cite[Theorem~1.1]{KoLa14}, $(X,d,\mu)$ supports a $(q,1)$-Poincar\'e inequality if and only if it supports a $(q,1)$-relative isoperimetric inequality, with quantitatively related constants.  Moreover, if $(X,d,\mu)$ is a doubling length space, supporting a $(q,1)$-Poincar\'e inequality, then as shown in \cite{HaKo00}, we may take $\tau=1$ to be the scaling constant in \eqref{eq:PI}, and consequently in \eqref{eq:rel isoperimetric}.  

Furthermore, if $(X,d,\mu)$ is a doubling metric measure space supporting a $(1,1)$-Poincar\'e inequality with scaling constant $\tau$, then by \cite[Theorem~4.21]{BB11} for example, $(X,d,\mu)$ also supports a $(\frac{Q}{Q-1},1)$-Poincar\'e inequality with scaling constant $2\tau$, where $Q:=Q_d$ is the exponent from \eqref{eq:rel lower mass bound exponent}.  Applying this inequality to Lipschitz approximating sequences, the following holds for all measurable sets $E\subset X$ and balls $B\subset X$:   
\begin{equation}\label{eq:Q rel isoperimetric ineq}
\left(\frac{\min\{\mu(B\cap E),\,\mu(B\setminus E)\}}{\mu(B)}\right)^{(Q-1)/Q}\le C\rad(B)\frac{ P
    (E,2\tau B)}{\mu(2\tau B)}.   
\end{equation}
Here, $C\ge 1$ depends only on the doubling constant and constants associated with the $(1,1)$-Poincar\'e inequality.

\subsection{Fractional perimeter and fractional inequalities}

Given $0<\theta<1$, a measurable set $E\subset X$, and an open set $\Omega\subset X$, we define the fractional perimeter of $E$ in $\Omega$ by 
\begin{equation}\label{eq:frac per}
    P_\theta(E,\Omega):=\int_{\Omega\cap E}\int_{\Omega\setminus E}\frac{2}{d(x,y)^\theta[\mu(B(x,d(x,y)))+\mu(B(y,d(x,y)))]}d\mu(y)d\mu(x).
\end{equation}
Note that this definition agrees, up to a dimensional normalization constant, with the Euclidean definition given by \eqref{eq:intro euc frac per}, and the symmetry of the kernel yields  $P_\theta(E,\Omega)=P_\theta(X\setminus E,\Omega)$.  If $\mu$ is doubling then it also follows that 
\[
P_\theta(E,\Omega)\simeq\int_{\Omega\cap E}\int_{\Omega\setminus E}\frac{1}{d(x,y)^\theta\mu(B(x,d(x,y)))}d\mu(y)d\mu(x).
\]

We have the following fractional coarea formula, see \cite{Vis91} and \cite[Lemma~4.3]{PS20} for such results in the Euclidean setting: 

\begin{lemma}\label{lem:frac coarea}
    Suppose that $\mu$ is $\sigma$-finite.  Let $u\in L^1_\loc(X)$, and let $\Omega\subset X$ be an open set. Then for each $0<\theta<1$, we have that 
    \[
    \int_\Omega\int_\Omega\frac{|u(x)-u(y)|}{d(x,y)^\theta[\mu(B(x,d(x,y)))+\mu(B(y,d(x,y)))]}d\mu(y)d\mu(x)=\int_\R P_\theta(\{u>t\},\Omega)dt.
    \]
\end{lemma}

\begin{proof}
    For $(\mu\times\mu)$-a.e.\ $(x,y)\in\Omega\times\Omega$, we have that 
    \[
    |u(x)-u(y)|=\int_\R|\chi_{\{u>t\}}(x)-\chi_{\{u>t\}}(y)|dt,
    \]
    and so by Tonelli's theorem and the definition \eqref{eq:frac per}, we have that 
    \begin{align*}
      \int_\Omega\int_\Omega&\frac{|u(x)-u(y)|}{d(x,y)^\theta[\mu(B(x,d(x,y)))+\mu(B(y,d(x,y)))]}d\mu(y)d\mu(x)\\
      &=\int_\R\int_\Omega\int_\Omega\frac{|\chi_{\{u>t\}}(x)-\chi_{\{u>t\}}(y)|}{d(x,y)^\theta[\mu(B(x,d(x,y)))+\mu(B(y,d(x,y)))]}d\mu(y)d\mu(x)dt\\
      &=\int_{\R}P_\theta(\{u>t\},\Omega)dt.\qedhere
    \end{align*}
\end{proof}

For the convenience of the reader, we conclude this section by listing all fractional Poincar\'e, isoperimetric, and Sobolev-type inequalities referenced in the introduction.  These will be the principal objects of study in what follows.

\begin{definition}\label{def:frac ineq}
Let $0<\theta<1$ and $1\le q<\infty$.
\begin{enumerate}
    \item We say that $(X,d,\mu)$ supports a \emph{$(\theta,q,1)_{BBM}$-Poincar\'e inequality} if there exist $C,\tau\ge 1$ such that for all balls $B:=B(x,r)\subset X$ and $u\in L^1_\loc(X)$, 
    \[
    \left(\fint_{B}|u-u_{B}|^qd\mu\right)^{1/q}\le C(1-\theta)r^\theta\fint_{\tau B}\int_{\tau B}\frac{|u(x)-u(y)|}{d(x,y)^\theta\mu(B(x,d(x,y)))}d\mu(y)d\mu(x).
    \]
     \item We say that $(X,d,\mu)$ supports a \emph{$(\theta,q,1)_{BBM}$-relative isoperimetric inequality} if there exist $C,\tau\ge 1$ such that for all balls $B:=B(x,r)\subset X$ and measurable $E\subset X$, 
    \[
    \left(\frac{\min\{\mu(B\cap E),\,\mu(B\setminus E)\}}{\mu(B)}\right)^{1/q}\le C(1-\theta)r^\theta\frac{ P_\theta(E,\tau B)}{\mu(\tau B)}.
    \]
    \item We say that $(X,d,\mu)$ supports an \emph{improved $(\theta,q,1)_{BBM}$-Poincar\'e inequality} if there exist $C,\tau\ge 1$ such that for all balls $B:=B(x,r)\subset X$ and $u\in L^1_\loc(X)$, 
    \[
    \left(\int_{B}|u-u_{B}|^qd\mu\right)^{1/q}\le C(1-\theta)\int_{\tau B}\int_{\tau B}\frac{|u(x)-u(y)|}{d(x,y)^\theta\mu(B(x,d(x,y)))}d\mu(y)d\mu(x).
    \]
     \item We say that $(X,d,\mu)$ supports an \emph{improved $(\theta,q,1)_{BBM}$-relative isoperimetric inequality} if there exist $C,\tau\ge 1$ such that for all balls $B:=B(x,r)\subset X$ and measurable $E\subset X$, 
    \[
    \min\{\mu(B\cap E),\,\mu(B\setminus E)\}^{1/q}\le C(1-\theta)P_\theta(E,\tau B).
    \]
    \item We say that $(X,d,\mu)$ supports a \emph{$(\theta,q,1)_{BBM}$-Sobolev inequality} if there exists $C\ge 1$ such that for all $u\in L^1(X)$, 
    \[    \left(\int_{X}|u|^qd\mu\right)^{1/q}\le C\theta(1-\theta)\int_{X}\int_{X}\frac{|u(x)-u(y)|}{d(x,y)^\theta\mu(B(x,d(x,y)))}d\mu(y)d\mu(x).
    \]
    \item We say that $(X,d,\mu)$ supports a \emph{$(\theta,q,1)_{BBM}$-isoperimetric inequality} if there exists $C\ge 1$ such that for all measurable $E\subset X$ with $\mu(E)<\infty$, 
    \[    \mu(E)^{1/q}\le C\theta(1-\theta)P_\theta(E,X).
    \]
\end{enumerate}
In each of these definitions, the constants $C$ and $\tau$ can depend on structural data associated to $X$, but are independent of $\theta$, $u$, $B$, and $E$. 
\end{definition}

\section{Nonlocal improvement of the relative isoperimetric inequality}\label{FRIIne}
The main result of this section, Lemma~\ref{lem:Fractional Isoperimetric}, is a nonlocal, isoperimetric-type inequality, which will be used to prove the fractional Poincar\'e inequality Theorem~\ref{thm: MFPinq} in the next section.  As discussed in Remark~\ref{improvedRiso}, the inequality given in Lemma~\ref{lem:Fractional Isoperimetric} acts as an improvement of the standard relative isoperimetric inequality \eqref{eq:Q rel isoperimetric ineq}.  The key lemma used in establishing this inequality is Lemma~\ref{lem:boundary balls}, which allows us to decompose a set near its boundary into a suitable collection of balls at a given, sufficiently small scale. 

 The lemmas of this section are metric setting analogs of \cite[Lemmas~3.1,\,3.2,\,3.3]{MPW24}, which were proven in the Euclidean setting and used to establish the corresponding fractional Poincar\'e inequalities there.  Our proofs of these lemmas are similar, but key differences appear in the argument for Lemma~3.1. In our setting, we do not have access to a system of dyadic cubes (whose perimeters are controlled), and so we must take special care to control the overlap of certain collections of balls.

As mentioned in the introduction, it suffices to prove Theorem~\ref{thm: MFPinq} under the additional assumption that $(X,d)$ is geodesic, see the discussion at the beginning of Section~\ref{FPI}. As such, we will also additionally assume that $(X,d)$ is geodesic in order to establish the lemmas in this section.  This assumption is quite useful, as it allows us to take the scaling constant of the $(1,1)$-Poincar\'e inequality, and hence $(1,1)$-relative isoperimetric inequality, to be $\tau=1$ \cite{HaKo00}, and it gives us access to the annular decay estimate Lemma~\ref{lem:annular decay}, for example.

\begin{lemma}\label{lem:boundary balls}
    Let $(X,d,\mu)$ be a complete, geodesic metric measure space, with $\mu$ a doubling measure supporting a $(1,1)$-Poincar\'e inequality. Let $B_0:=B(x_0,r_0)\subset X$, and let $E\subset X$ be a measurable set such that 
    \[
    \lambda\le\frac{\mu(B_0\cap E)}{\mu(B_0)}\le 1-\lambda
    \]
    for some $0<\lambda\le 1/2$.
    Then there exist constants $C:=C(C_\mu)\ge 1$ and $K:=K(\lambda,C_\mu)\in\N$ so that for all $k\in\N$, $k\ge K$, there exists a countable collection of balls $\{B_i\}_{i\in\N}$, with $B_i\subset B_0$,
    $2^{-k}r_0\le \rad(B_i)\leq 2^{-k+4}r_0$, and $\{\frac{1}{60}B_i\}_{i\in\N}$ pairwise disjoint, such that 
    \begin{align}
        &(i)\quad \lambda/C\le\frac{\mu(B_i\cap E)}{\mu(B_i)}\le 1-\lambda/C,\label{eq:density}\\
       &(ii)\quad 2^{-k}\mu(B_0)\lesssim\frac{1}{\lambda}\sum_i\mu(B_i),\label{eq:total mass}
    \end{align}
 with the comparison constant in $(ii)$ depending only on $C_\mu$ and the constants associated with the $1$-Poincar\'e inequality.  We may take $C=16C_\mu^{19}$ and 
 \[
 K=\left\lceil\log_2\left(32\left(4C_A/\lambda\right)^{1/\beta}\right)\right\rceil,
 \]  
 where $C_A\ge 1$ and $0<\beta\le 1$ are the constants from Lemma~\ref{lem:annular decay}, depending only on $C_\mu$.
\end{lemma}

\begin{proof}
   Let $k\in\N$ be such that $k\ge K$, with $K\in\N$ as defined above. Since $\mu$ is doubling, there exists a countable collection of balls $\{B(x_i,2^{-k}r_0)\}_i$ covering $X$ such that $\{B(x_i,2^{-k}r_0/5)\}_i$ is pairwise disjoint.  For each $i$, there exists $2^{-k}r_0\le \rho_i<2^{-k+1}r_0$, by \cite[Lemma~6.2]{KKST}, such that 
    \begin{equation}\label{eq:ball perimeter}    P(B(x_i,\rho_i),X)\lesssim\frac{\mu(B(x_i,\rho_i))}{2^{-k}r_0}.
    \end{equation}
    Let $I:=\{i\in\N:B(x_i,\rho_i)\subset B_0\}$ and denote $\B_k:=\{B(x_i,\rho_i)\}_{i\in I}$.  We note that the collection $\{\frac{1}{10}B\}_{B\in\B_k}$ is pairwise disjoint.

    Let $C_0:=16C_\mu^{16}$, and define the subcollections
    \begin{equation*}
    \mathcal D_k^1:=\left\{B\in\B_k:\frac{\mu(B\cap E)}{\mu(B)}\ge \lambda/C_0\right\}
    \end{equation*}
	and
	\begin{equation*}
		\mathcal D_k^2:=\left\{B\in\B_k:\frac{\mu(B\setminus E)}{\mu(B)}\ge \lambda/C_0\right\}.
	\end{equation*}
    For ease of notation, we set
    \[
    A_1:=\bigcup_{B\in\Dd_k^1}B,\quad A_2:=\bigcup_{B\in\Dd_k^2}B,\quad A_1':=\bigcup_{B\in\Dd_k^1\setminus\Dd_k^2}B,\quad A_2':=\bigcup_{B\in\Dd_k^2\setminus\Dd_k^1}B.
    \]
    By the definition of ${\mathcal{B}_k}$, we have that $B(x_0,(1-2^{-k+2})r_0)\subset\bigcup_{B\in\B_k}B$. 
    Thus, by Lemma~\ref{lem:annular decay} and by our choice of $K$ and $C_0$, we have that
    \begin{align*}
    \mu\left(E\cap B_0\setminus A_1\right)
    &\le\mu\left(B_0\setminus\bigcup_{B\in\B_k}B\right)+\sum_{B\in \mathcal B_k\setminus \mathcal D_k^1}\mu(E\cap B)\\
    &\le \mu\left(B(x_0,r_0)\setminus B(x_0,(1-2^{-k+2})r_0)\right)+\lambda C_0^{-1}\sum_{B\in \mathcal B_k\setminus \mathcal D_k^1}\mu(B)\\
    &\le C_A (2^{-k+2})^\beta\mu(B_0)+\lambda C_\mu^4 C_0^{-1}\sum_{B\in\B_k}\mu\left(\frac{1}{10}B\right)\\
    &\le (C_A(2^{-k+2})^\beta+\lambda C_\mu^4C_0^{-1})\mu(B_0)\\
    &\le\frac{\lambda}{2}\mu(B_0).
    \end{align*}
	We also have for the complement $E^c$
	  \begin{align*}
		\mu\left(E^c\cap B_0\setminus A_1\right)
		\le \mu\left(E^c\cap B_0\right)
		\le \left(1-\lambda\right)\mu(B_0),
	\end{align*}
    and so combining these two estimates, we obtain
    \[
    \mu\left(B_0\setminus A_1\right)\le \left(1-\lambda/2\right)\mu(B_0).
    \]
    Thus
	\begin{equation}\label{eq:A_1 measure}
	\mu\left(A_1\right)
	\ge\frac{\lambda}{2}\mu(B_0).
	\end{equation}
	By an analogous argument, replacing $\Dd_k^1$ with $\Dd_k^2$, we also have
	\begin{equation}\label{eq:A_2 measure}
	\mu\left(A_2\right)
	\ge \frac{\lambda}{2}\mu(B_0).
	\end{equation}
    
    In the case that
    \[
	\mu\left(\bigcup_{B\in \mathcal D_k^1\cap \mathcal D_k^2}B\right)
	\ge\frac{\lambda}{8}\mu(B_0),
	\]
	the collection of balls $\Dd_k^1\cap\Dd_k^2$ satisfies the conclusion of the lemma with $C=16C_\mu^{16}$, by the definition of $\Dd_k^1$ and $\Dd_k^2$.  We therefore assume that 
    \[
	\mu\left(\bigcup_{B\in \mathcal D_k^1\cap \mathcal D_k^2}B\right)
	<\frac{\lambda}{8}\mu(B_0).
    \]
    \vskip.2cm
    
    \noindent{\it Claim:} The following holds for either $i=1$ or $i=2$:
	\begin{align*}
	\mu\left(A_i\right)
	\le \left(1/2+\lambda/4\right)\mu(B_0).
    \end{align*}
	
    \vskip.2cm
    \noindent{\it Proof of claim:} Suppose that this inequality fails for both $A_1$ and $A_2$. We then have that

    \begin{align}\label{eq:Intersection estimate}
        (1+\lambda/2)\mu(B_0)\le\mu(A_1)+\mu(A_2)&\le\mu(A_1')+\mu(A_2')+2\mu\left(\bigcup_{B\in \mathcal D_k^1\cap\Dd_k^2}B\right)\nonumber\\
        &\le\mu(A_1'\cup A_2')+\mu(A_1'\cap A_2')+\frac{\lambda}{4}\mu(B_0)\nonumber\\
        &\le (1+\lambda/4)\mu(B_0)+\mu(A_1'\cap A_2').
    \end{align}
    If $B\in\Dd_k^1\setminus\Dd_k^2$ and $B'\in\Dd_k^2\setminus\Dd_k^1$, then it follows that 
    \[
    \frac{\mu(B\cap E)}{\mu(B)}\ge(1-\lambda/C_0)\quad\text{and}\quad\frac{\mu(B'\cap E)}{\mu(B')}<\lambda/C_0.
    \]
    If $B\cap B'\ne\varnothing$, then $B'\subset 5B$, and so it follows from doubling that 
    \begin{align*}
        \frac{\mu(B\cap B')}{\mu(B)}&=1-\frac{\mu(B\cap E)}{\mu(B)}+\frac{\mu(B\cap E\cap B')}{\mu(B)}-\frac{\mu(B\setminus(E\cup B'))}{\mu(B)}\\
        &\le\lambda/C_0+\frac{\mu(B'\cap E)}{\mu(B)}\le\lambda/C_0+ C_\mu^3\frac{\mu(B'\cap E)}{\mu(B')}\le 2\lambda C_\mu^3/C_0.
    \end{align*}
   Therefore, we have that 
    \begin{align*}
        \mu(A_1'\cap A_2')\le\sum_{B\in\Dd_k^1\setminus\Dd_k^2}\mu(B\cap A_2')\le&\sum_{B\in\Dd_k^1\setminus\Dd_k^2}\sum_{B'\in\Dd_k^2\setminus\Dd_k^1}\mu(B\cap B')
        \\
        &\le 2\lambda C_\mu^3C_0^{-1} \sum_{B\in\Dd_k^1\setminus\Dd_k^2}\sum_{\substack{ B'\in\Dd_k^2\setminus\Dd_k^1\\ B'\cap B\ne\varnothing}}\mu(B).
    \end{align*}
    Since $\mu$ is doubling and the collection $\{\frac{1}{10}B\}_{B\in\B_k}$ is pairwise disjoint, it follows that for each $B\in\Dd_k^1\setminus\Dd_k^2$, there are at most $C_\mu^9$ balls in $\Dd_k^2\setminus\Dd_k^1$ which intersect $B$.  Thus, by our choice of $C_0$, we have that 
    \begin{align*}
        \mu(A_1'\cap A_2')\le 2\lambda C_\mu^{12}C_0^{-1}\sum_{B\in\Dd_k^1\setminus\Dd_k^2}\mu(B)\le 2\lambda C_\mu^{16}C_0^{-1}\sum_{B\in\Dd_k^1\setminus\Dd_k^2}\mu\left(\frac{1}{10}B\right)&\le 2\lambda C_\mu^{16}C_0^{-1}\mu(B_0)\\
        &\le\frac{\lambda}{8}\mu(B_0).
    \end{align*}
    Combining this inequality with \eqref{eq:Intersection estimate} yields a contradiction, proving the claim.
    \vskip.2cm

    We now assume that the claim holds with $i=1$.  In this case, we define the collection
    \[
    \CC_k:=\{B\in\Dd_k^1:\text{ there exists }B'\in\B_k\setminus\Dd_k^1\st 2B\cap B'\ne\varnothing,\,6B\subset B_0\},
    \]
    and we show that the collection $\{6B\}_{B\in\CC_k}$ satisfies the conclusion of the lemma.

    To this end, it follows from the definition of $\CC_k$ that for all $B\in\CC_k$, there exists $B'\in\B_k\setminus\Dd_k^1$ such that $B'\subset 6B$.  Thus, by doubling, and the fact that $6B\subset 16B'$, we have that
    \begin{align}\label{eq:Case 1 density}
        \frac{\lambda}{C_\mu^3C_0}\le\frac{\mu(B\cap E)}{C_\mu^3\mu(B)}\le\frac{\mu(6B\cap E)}{\mu(6B)}&\le\frac{\mu(6B\setminus B')}{\mu(6B)}+\frac{\mu(B'\cap E)}{\mu(B')}\\
        &\le\left(1-\frac{\mu(B')}{\mu(6B)}\right)+\frac{\lambda}{C_0}\nonumber\le 1-\frac{1}{C_\mu^4}+\frac{\lambda}{C_0}\le 1-\frac{1}{2C_\mu^4},
    \end{align}
    where the last inequality follows from our choice of $C_0$.  Hence, $\{6B\}_{B\in\CC_k}$ satisfies $(i)$ with $C=16C_\mu^{19}$.
    
By the relative isoperimetric inequality \eqref{eq:rel isoperimetric}, it follows that 
\begin{align*}
    \min\{\mu((1-2^{-k+5})B_0\cap A_1),\,\mu((1-2^{-k+5})B_0\setminus A_1)\}\lesssim r_0 P\left(A_1,(1-2^{-k+5})B_0\right).
\end{align*}
Note that as $(X,d)$ is geodesic, the scaling constant of the relative isoperimetric inequality \eqref{eq:rel isoperimetric} is taken here to be $\tau=1$. If $x\in(\partial A_1)\cap(1-2^{-k+5})B_0$, then there exists $B\in\Dd_k^1$ such that $x\in\partial B$ and $6B\subset B_0$. Moreover, since $(1-2^{-k+2})B_0\subset\bigcup_{B\in\B_k}B$, there exists $B'\in\B_k\setminus\Dd_k^1$ such that $x\in B'$.  Therefore, $B\in\CC_k$, and so we have that 
\[
(\partial A_1)\cap(1-2^{-k+5})B_0\subset\bigcup_{B\in\CC_k}\partial B.
\]
From this fact, as well as \eqref{eq:ball perimeter}, it follows from the previous inequality that 
\begin{align}\label{eq:A_1 isoperimetric}
    \min\{\mu((1-2^{-k+5})B_0\cap A_1),\,\mu((1-2^{-k+5})B_0\setminus A_1)\}
    &\lesssim r_0 \sum_{B\in\CC_k} P(B,X)\nonumber\\
    &\lesssim r_0\sum_{B\in\CC_k}\frac{\mu(B)}{2^{-k}r_0}\nonumber\\
    &\le 2^k\sum_{B\in\CC_k}\mu(6B). 
\end{align}
By our choice of $K$ and \eqref{eq:A_1 measure}, we have that
\begin{align*}
    \mu((1-2^{-k+5})B_0\cap A_1)\ge\mu(A_1)-\mu(B_0\setminus(1-2^{-k+5})B_0)&\ge\frac{\lambda}{2}\mu(B_0)-C_A(2^{-k+5})^\beta\mu(B_0)\\
    &\ge\frac{\lambda}{4}\mu(B_0).
\end{align*}
Likewise, by the claim and our choice of $K$, we have that 
\begin{align*}
    \mu((1-2^{-k+5})B_0\setminus A_1)&\ge\mu(B_0\setminus A_1)-\mu(B_0\setminus(1-2^{-k+5})B_0)\\
    &\ge(1/2-\lambda/4)\mu(B_0)-C_A(2^{-k+5})^\beta\mu(B_0)\ge\frac{1}{4}\mu(B_0).
\end{align*}
Therefore, by \eqref{eq:A_1 isoperimetric}, we have that
\begin{align*}
    2^{-k}\mu(B_0)\lesssim\frac{1}{\lambda}\sum_{B\in\CC_k}\mu(6B),
\end{align*}
and so $\{6B\}_{B\in\CC_k}$ satisfies $(ii)$.

In the case that the claim holds with $i=2$, we define the collection $\CC_k$ to be 
 \[
    \CC_k:=\{B\in\Dd_k^2:\text{ there exists }B'\in\B_k\setminus\Dd_k^2\st 2B\cap B'\ne\varnothing,\,6B\subset B_0\},
    \]
and similarly show that the collection $\{6B\}_{B\in\CC_k}$ satisfies the conclusion of the lemma.  Indeed, by the same argument used to obtain \eqref{eq:Case 1 density} in the previous case, as well as our choice of $C_0$, it follows that for all $B\in\CC_k$, 
\begin{align*}
    \frac{1}{2C_\mu^4}\le\frac{1}{C_\mu^{4}}-\frac{\lambda}{C_0}\le\frac{\mu(6B\cap E)}{\mu(6B)}\le1-\frac{\lambda}{C_\mu^3C_0},
\end{align*}
and so $\{6B\}_{B\in\CC_k}$ satisfies $(i)$ with $C=16C_\mu^{19}$.  Likewise, by \eqref{eq:A_2 measure}, the claim, and by applying the relative isoperimetric inequality as above with respect to $A_2$ rather than $A_1$, we obtain
\begin{align*}
    2^{-k}\mu(B_0)\lesssim\frac{1}{\lambda}\sum_{B\in\CC_k}\mu(6B),
\end{align*}
again using our choice of $K$.  Thus, $(ii)$ is satisfied in this case, concluding the proof.
\end{proof}

\begin{lemma}\label{lem:Annuli estimate}
	Let $(X,d,\mu)$ be a geodesic metric measure space with $\mu$ a doubling measure. Let $0<r_0<2\diam(X)$, $B_0:=B(x_0,r_0)\subset X$, $a \leq r_0 / 2$, and $0<\varepsilon\le1/2$. Let $B_1:=B(x_1,r_1)\subset B_0$ with 
	\begin{equation*}
	r_1 \le\frac{a}{2(32C_\mu^4C_A)^{1/\beta}}.
	\end{equation*}
    Then for any measurable set $E\subset X$ with 
	\begin{equation}\label{eq:eps density}
	\varepsilon \leq 
	\frac{\mu(B_1\cap  E)}{\mu(B_1)}
	\leq 1-\varepsilon,
	\end{equation}
	we have 
	\[
	\mu(B_1)\leq 
	\frac{4}{\varepsilon}\int_{B_1}\left\lvert\chi_{E}(x)-\frac{\mu(A(x)\cap E)}{\mu(A(x))} \right\rvert d\mu(x),
	\]
	where $A(x)=B_0 \cap B(x,a)\backslash B(x,a/2)$.
\end{lemma}

\begin{proof}
    Since $0<r_0<2\diam(X)$ and since $X$ is geodesic, it follows that $A(x)\ne\varnothing$ for each $x\in B_1$.  We first claim that for each $x\in B_1$, we have
    \begin{align}\label{eq:1/4 estimate}
        \left|\frac{\mu(A(x)\cap E)}{\mu(A(x))}-\frac{\mu(A(x_1)\cap E)}{\mu(A(x_1))}\right|\le\frac{1}{4}.
    \end{align}
    To prove this, for each $x\in B_1$, set 
    \[
    \wtil A(x):= B(x,a)\setminus B(x,a/2).
    \]
    We then have that
    \begin{align*}
        \left|\mu(A(x)\cap E)-\mu(A(x_1)\cap E)\right|&=|\mu(A(x)\cap E\setminus A(x_1))-\mu(A(x_1)\cap E\setminus A(x))|\\
        &\le\mu(A(x)\setminus A(x_1))+\mu(A(x_1)\setminus A(x))\\
        &\le\mu(\wtil A(x)\setminus\wtil A(x_1))+\mu(\wtil A(x_1)\setminus\wtil A(x)).
    \end{align*}
    Since 
    \[
    \wtil A(x)\setminus\wtil A(x_1)\subset\Big(B(x,a/2+d(x,x_1))\setminus B(x,a/2)\Big)\cup\Big(B(x,a)\setminus B(x,a-d(x,x_1))\Big),
    \]
    it follows from Lemma~\ref{lem:annular decay} and our assumption on $r_1$ that 
    \begin{align*}
    \mu(\wtil A(x)&\setminus\wtil A(x_1))\\
    &\le C_A\left(\frac{d(x,x_1)}{a/2+d(x,x_1)}\right)^\beta\mu(B(x,a/2+d(x,x_1)))+C_A\left(\frac{d(x,x_1)}{a}\right)^\beta\mu(B(x,a))\\
    &\le 2C_A\left(\frac{2r_1}{a}\right)^\beta\mu(B(x,a)).
    \end{align*}
    By a similar argument and doubling, it also follows that 
    \begin{align*}
        \mu(\wtil A(x_1)\setminus\wtil A(x))\le 2C_A\left(\frac{2r_1}{a}\right)^\beta\mu(B(x_1,a))\le2C_\mu C_A\left(\frac{2r_1}{a}\right)^\beta\mu(B(x,a)).
    \end{align*}
    Therefore, we have that 
    \begin{align*}
         \left|\mu(A(x)\cap E)-\mu(A(x_1)\cap E)\right|&\le\mu(A(x)\setminus A(x_1))+\mu(A(x_1)\setminus A(x))\\
         &\le 4C_\mu C_A\left(\frac{2r_1}{a}\right)^\beta\mu(B(x,a)).
    \end{align*}
    From this inequality, we then obtain
    \begin{align*}
        |\mu(A(x)&\cap E)\mu(A(x_1))-\mu(A(x_1)\cap E)\mu(A(x))|\\
        &\le|\mu(A(x)\cap E)\mu(A(x_1))-\mu(A(x_1)\cap E)\mu(A(x_1))|\\
        &\qquad+|\mu(A(x_1)\cap E)\mu(A(x_1))-\mu(A(x_1)\cap E)\mu(A(x))|\\
        &=\left|\mu(A(x)\cap E)-\mu(A(x_1)\cap E)\right|\mu(A(x_1))+|\mu(A(x_1))-\mu(A(x))|\mu(A(x_1)\cap E)\\
        &\le 8C_\mu C_A\left(\frac{2r_1}{a}\right)^{\beta}\mu(B(x,a))\mu(A(x_1)).
    \end{align*}

We note that there exists $y\in A(x)$ so that $B(y,a/4)\subset A(x)$.  Indeed, if $x\in B(x_0,r_0/2)$, then this follows by connectedness of $X$ and since $B(x,a)\subset B_0$.  If $x\in B_0\setminus B(x_0,r_0/2)$, then we join $x$ to $x_0$ by a geodesic, and choose $y$ to be a point on this geodesic such that $d(x,y)=3a/4$.  Therefore, in either case, we have by doubling that 
\begin{align*}
\mu(B(x,a))\le \mu(B(y,2a))\le C_\mu^3\mu(B(y,a/4))\le C_\mu^3\mu(A(x)).
\end{align*}
Combining this with the previous inequality, we then obtain
\begin{align*}
    |\mu(A(x)&\cap E)\mu(A(x_1))-\mu(A(x_1)\cap E)\mu(A(x))|\le 8C_\mu^4 C_A\left(\frac{2r_1}{a}\right)^{\beta}\mu(A(x))\mu(A(x_1)).
\end{align*}
Our assumption on $r_1$ then gives us \eqref{eq:1/4 estimate}, proving the claim.

If $\mu(A(x_1)\cap E)/\mu(A(x_1))\ge 1/2$, then \eqref{eq:1/4 estimate} implies that $\mu(A(x)\cap E)/\mu(A(x))\ge 1/4$.  In this case, we have from \eqref{eq:eps density} that
\begin{align*}
    \mu(B_1)\le\frac{1}{\eps}\mu(B_1\setminus E)=\frac{1}{\eps}\int_{B_1\setminus E}1\,d\mu\le\frac{4}{\eps}\int_{B_1\setminus E}\left|\chi_E(x)-\frac{\mu(A(x)\cap E)}{\mu(A(x))}\right|d\mu(x).
\end{align*}
Likewise, if $\mu(A(x_1)\cap E)/\mu(A(x_1))<1/2$, then \eqref{eq:1/4 estimate} implies that $1-\mu(A(x)\cap E)/\mu(A(x))\ge 1/4$, in which case \eqref{eq:eps density} yields
\begin{align*}
    \mu(B_1)\le\frac{1}{\eps}\mu(B_1\cap E)=\frac{1}{\eps}\int_{B_1\cap E}1\,d\mu\le\frac{4}{\eps}\int_{B_1\cap E}\left|\chi_E(x)-\frac{\mu(A(x)\cap E)}{\mu(A(x))}\right|d\mu(x).
\end{align*}
Therefore, the conclusion of the lemma holds in either case.
\end{proof}

\begin{lemma}\label{lem:Fractional Isoperimetric}
Let $(X,d,\mu)$ be a complete geodesic metric measure space, with $\mu$ a doubling measure supporting a $(1,1)$-Poincar\'e inequality.  Let $0<r<2\diam(X)$, $B:=B(x_0,r) \subset X$ be a ball, $E\subset X$ a measurable set, and $k\in\mathbb{N}$ be such that
\[
\frac{1}{2^{kQ}}\leq\frac{\mu(B\cap E)}{\mu(B)}\leq \frac{1}{2},
\]
where $Q:=Q_d>1$ is the relative lower mass bound exponent \eqref{eq:rel lower mass bound exponent}.
Then there exists $C\ge 1$ such that
\begin{equation}\label{eq:FracIso}
\left(\frac{\mu(B\cap E)}{\mu(B)}\right)^{\frac{Q-1}{Q}}
\leq 
C 2^{k} \fint_{B}\fint_{B\cap B(x,2^{-k}r)\backslash B(x,2^{-k-1}r)}|\chi_E(x)-\chi_{E}(y)|d\mu(y)d\mu(x).
\end{equation}
Here $C$ depends only on the doubling and $(1,1)$-Poincar\'e inequality constants.
\end{lemma}

\begin{proof}
Since $\mu(B\cap E)/\mu(B)\le 1/2$, we apply Lemma~\ref{lem:the decomposition of the ball} to obtain a countable collection of pairwise disjoint balls $B_i\subset B$ such that 
\begin{align*}
    &(i)\quad r_i:=\rad(B_i)=2^{-M_i}r\text{ for some }M_i\in\N\cup\{0\},\\
	&(ii)\quad B \cap E \subset \bigcup_{i}5B_i \quad \textrm{up to a set of }\mu\textrm{-measure zero},\\	
	&(iii)\quad \frac{1}{2^{k_0Q+1}}\le\frac{\mu(B_i\cap E)}{\mu(B_i)}\leq \frac{1}{2},
	\end{align*}
where $k_0$ is the smallest positive integer with $2^{-k_0Q}\leq C_\mu^{-2}$.  Denote by $k_1$ the smallest positive integer such that 
\begin{align}\label{eq:k_1}
2^{-k_1}\le\frac{1}{32(32 C_\mu^4C_A)^{1/\beta}},
\end{align}
and let 
\[
K:=\left\lceil\log_2\left(32(4C_A2^{k_0Q+1})^{1/\beta}\right)\right\rceil.
\]
We apply Lemma~\ref{lem:boundary balls} with $\max\{k+k_1-M_i,K\}\in\N$ for $E$ on each $B_i$ to obtain a collection $\{B_{i,j}\}_{j=1}^{N_i}$ of balls, with $B_{i,j}\subset B_i$ and $\{\frac{1}{60}B_{i,j}\}_{j=1}^{N_i}$ pairwise disjoint, such that  
\begin{align}\label{eq:r_{i,j}}
2^{-\max\{k+k_1-M_i,K\}}r_i\le r_{i,j}:=\rad(B_{i,j})\le 2^{-\max\{k+k_1-M_i,K \}+4}r_i,
\end{align}
and such that
\begin{align}\label{eq:Bij density}
\frac{1}{2^{k_0Q+1}(16C_\mu^{19})}\leq\frac{\mu(B_{i,j}\cap E)}{\mu(B_{i,j})}\leq 1-\frac{1}{2^{k_0Q+1}(16C_\mu^{19})}.
\end{align}
From this lemma, we also have that 
\[
2^{-\max\{k+k_1-M_i, K\}}\mu(B_i)\lesssim\sum_j\mu(B_{i,j}).
\]
Using these properties of $B_i$ and $B_{i,j}$, as well as doubling, it then follows that
\begin{align}\label{eq:double sum density}
\left(\frac{\mu(B\cap E)}{\mu(B)}\right)^{\frac{Q-1}{Q}}
&\lesssim\left(\frac{\mu(B\cap E)}{\mu(B)}\right) ^{-\frac{1}{Q}}\mu(B)^{-1}\sum_i \mu(B_i)\nonumber\\
&=\mu(B)^{-1}\min\left\{\left(\frac{\mu(B\cap E)}{\mu(B)}\right) ^{-\frac{1}{Q}},\,2^{k}\right\} \sum_i\mu(B_i)\nonumber\\
&\leq \mu(B)^{-1}\sum_i\min\left\{\left(\frac{\mu(B_i\cap E)}{\mu(B)}\right) ^{-\frac{1}{Q}},\,2^{k}\right\} \mu(B_i)\nonumber\\
&\leq \mu(B)^{-1}\sum_i\min\left\{2^{k_0+1/Q}\left(\frac{\mu(B_i)}{\mu(B)}\right) ^{-\frac{1}{Q}},\,2^{k}\right\} \mu(B_i)\nonumber\\
&\lesssim\mu(B)^{-1} \sum_i\min\left\{2^{k_0+1/Q+M_i},\,2^{k}\right\} \mu(B_i)\nonumber\\
&= 2^{k_0+1/Q+k+k_1+K}\mu(B)^{-1}\sum_i \min\left\{2^{M_i-k-k_1-K},\,2^{-k_0-1/Q-k_1-K}\right\} \mu(B_i)\nonumber\\
&\leq 2^{k_0+1/Q+k+k_1+K}\mu(B)^{-1}\sum_i 2^{-\max\{k+k_1-M_i,\,K\}}\mu(B_i)\nonumber\\
&\lesssim2^{k}\mu(B)^{-1}\sum_i\sum_j\mu(B_{i,j}),
\end{align}
where in the last inequality, we recall that $k_0$, $k_1$, and $K$ depend only on $C_\mu$.

Letting $a=2^{-k}r$, it follows from \eqref{eq:k_1} and \eqref{eq:r_{i,j}} that
\[
\rad(B_{i,j})\le\frac{a}{2(32C_\mu^4C_A)^{1/\beta}}.
\]
Thus, from \eqref{eq:Bij density}, we apply Lemma~\ref{lem:Annuli estimate} on $B_{i,j}$, with $\eps=(2^{k_0Q+1}(16C_\mu^{19}))^{-1}$, to obtain
\begin{align*}
    \mu(B_{i,j})\lesssim\int_{B_{i,j}}\left\lvert\chi_{E}(x)-\frac{\mu(A(x)\cap E)}{\mu(A(x))} \right\rvert d\mu(x),
\end{align*}
 where $A(x):=B\cap B(x,2^{-k}r)\backslash B(x,2^{-k-1}r)$. We then obtain
\begin{align*}
\sum_i\sum_j\mu(B_{i,j})
&\lesssim\sum_i\sum_j\int_{B_{i,j}}\left\lvert\chi_{E}(x)-\frac{\mu(A(x)\cap E)}{\mu(A(x))} \right\rvert d\mu(x)\nonumber\\
&=\sum_i\sum_j\int_{B_{i,j}}\left\lvert\chi_{E}(x)-\fint_{A(x)}\chi_{E}(y)d\mu(y) \right\rvert d\mu(x)\nonumber\\
&\leq \sum_i\sum_j\int_{B_{i,j}}\fint_{A(x)}\left\lvert \chi_E(x)-\chi_{E}(y)\right\rvert d\mu(y)d\mu(x) \nonumber\\
&\lesssim\sum_i\int_{B_i}\fint_{A(x)}\left\lvert \chi_E(x)-\chi_{E}(y)\right\rvert d\mu(y)d\mu(x)\nonumber\\
&\le\int_{B}\fint_{A(x)}|\chi_E(x)-\chi_E(y)|d\mu(y)d\mu(x).
\end{align*}
Here, the second to last inequality follows since the collection $\{B_{i,j}\}_{j=1}^{N_i}$ has bounded overlap, as $\{\frac{1}{60}B_{i,j}\}_{j=1}^{N_i}$ is pairwise disjoint.  The last inequality follows from pairwise disjointness of $\{B_i\}_i$.  
Combining this estimate with \eqref{eq:double sum density} yields the desired inequality.
\end{proof}

\begin{remark}\label{improvedRiso}
    In \cite[Remark~3.4]{MPW24}, the authors showed that in the Euclidean setting, the inequality \eqref{eq:FracIso} acts as a nonlocal improvement of the standard relative isoperimetric inequality \eqref{eq:Q rel isoperimetric ineq}. In this remark, we adapt their argument to show that the same holds true in the metric setting under the assumptions of the previous lemma. In particular, we show that for all measurable $E\subset X$, $0<r<2\diam(X)$, and $B:=B(x_0,r)\subset X$, the following holds for sufficiently large $k\in\N$:  
    \begin{align}\label{eq:isoper intermediary}
        \Bigg(\frac{\min\{\mu(B\cap E),\,\mu(B\setminus E)\}}{\mu(B)}\Bigg)^{(Q-1)/Q}
        &\le C2^k\fint_B\fint_{ A_k(x)}|\chi_E(x)-\chi_E(y)|d\mu(y)d\mu(x)\nonumber\\
        &\le C r\frac{P(E,(1+2^{-k+2})B)}{\mu(B)},
    \end{align}
    where $A_k(x)=B\cap B(x,2^{-k}r)\setminus B(x,2^{-k-1}r)$, and
    with $C\ge 1$ depending only on the doubling constant and constants associated with the $(1,1)$-Poincar\'e inequality.
    
    To prove this, we first note that 
    \[
    \fint_B\fint_{A_k(x)}|\chi_E(x)-\chi_E(y)|d\mu(y)d\mu(x)=\fint_B\fint_{A_k(x)}|\chi_{E^c}(x)-\chi_{E^c}(y)|d\mu(y)d\mu(x)
    \]
    and 
    \[
    P(E,(1+2^{-k+2})B)=P(E^c,(1+2^{-k+2})B).
    \]
    Furthermore, if either $\mu(B\cap E)=0$ or $\mu(B\setminus E)=0$, then \eqref{eq:isoper intermediary} holds trivially. Therefore, replacing $E$ with its complement $E^c$ if necessary, we may assume without loss of generality that 
    \[
    0<\frac{\mu(B\cap E)}{\mu(B)}\le\frac{1}{2}.
    \]
    Thus, there exists $K\in\N$ so that for all $k\in\N$ with $k\ge K$, we have that 
    \[
    \frac{1}{2^{kQ}}\le\frac{\mu(B\cap E)}{\mu(B)}\le\frac{1}{2}.
    \]
    Applying Lemma~\ref{lem:Fractional Isoperimetric} then yields the first inequality of \eqref{eq:isoper intermediary}.

To prove the second inequality in \eqref{eq:isoper intermediary}, we fix $k\ge K$ and cover $B$ with a collection $\{B_i\}_{i=1}^{N_k}$ of balls centered in $B$ such that $\rad(B_i)=2^{-k}r$ and such that $\{\frac{1}{5}B_i\}_{i=1}^{N_k}$ is pairwise disjoint.  We then have that 
\begin{align*}
\int_B\fint_{A_k(x)}|\chi_E(x)-&\chi_E(y)|d\mu(y)d\mu(x)\le\sum_{i=1}^{N_k}\int_{B_i}\fint_{A_k(x)}|\chi_E(x)-\chi_E(y)|d\mu(y)d\mu(x)\\
&\lesssim\sum_{i=1}^{N_k}\frac{1}{\mu(B_i)}\left(\int_{B_i\cap E}\mu(A_k(x)\setminus E)d\mu(x)+\int_{B_i\setminus E}\mu(A_k(x)\cap E)d\mu(x)\right)\\
&\lesssim\sum_{i=1}^{N_k}\frac{\mu(2B_i\cap E)\mu(2B_i\setminus E)}{\mu(B_i)}\\
&\lesssim\sum_{i=1}^{N_k}\min\{\mu(2B_i\cap E),\,\mu(2B_i\setminus E)\}.
\end{align*}
Here we have used the assumption that $X$ is geodesic and doubling, in particular to ensure that $\mu(A_k(x))\simeq\mu(B(x,2^{-k}r))\simeq\mu(B_i)$ for each $x\in B_i$.  See the proof of \eqref{eq:annulus sub ball} for the precise argument of this last statement.  Applying the relative isoperimetric inequality \eqref{eq:rel isoperimetric} and using the fact that $\{B_i\}_{i=1}^{N_k}$ has bounded overlap, it then follows that 
\begin{align*}
\int_B\fint_{A_k(x)}|\chi_E(x)-\chi_E(y)|d\mu(y)d\mu(x)\lesssim 2^{-k+1}r\sum_{i=1}^{N_k} P(E,2B_i)\lesssim 2^{-k}r P(E,(1+2^{-k+1})B),
\end{align*}
which gives us the second inequality in \eqref{eq:isoper intermediary}.
\end{remark}

\section{Fractional Poincar\'e inequality}\label{FPI} 

In this section, we prove Theorem~\ref{thm: MFPinq}, Corollary~\ref{thm: MFRiso}, and the equivalence given by Theorem~\ref{rmk: PiRiequiv}. Our proof of Theorem~\ref{thm: MFPinq} is inspired by the proof of \cite[Theorem~4.1]{MPW24}, established in the Euclidean setting.  In particular, a key step in the argument is the inequality given by Lemma~\ref{lem:Fractional Isoperimetric}.

As mentioned in the introduction, it suffices to prove Theorem~\ref{thm: MFPinq} under the additional assumption that $(X,d)$ is geodesic; we restate Theorem~\ref{thm: MFPinq} with this assumption below in Theorem~\ref{thm:MFPinq-geodesic}.  Indeed, if $(X,d,\mu)$ is a complete metric measure space, with $\mu$ a doubling measure supporting a $(1,1)$-Poincar\'e inequality, then $(X,d)$ is $L$-bilipschitz equivalent to $(X,\til d)$, where $\til d$ is the length metric induced from $d$, and $L$ depends only on the doubling and $(1,1)$-Poincar\'e inequality constants, see for example \cite[Corollary~8.3.16]{HKST}.  Moreover, completeness, the doubling property of $\mu$, and support of a $(1,1)$-Poincar\'e inequality are all preserved under a bilipschitz change of metric, with constants depending only on the original data, and consequently $(X,\til d)$ is geodesic.  Note that the relative lower mass bound exponent $Q_d$ given by \eqref{eq:rel lower mass bound exponent} is also preserved under a bilipschitz change of metric. Hence, given a ball $B_d(x,r)\subset X$ (with respect to $d$), $u\in L^1_\loc(X)$, and $1\le q\le Q_d/(Q_d-\theta)$, it follows from doubling and Theorem~\ref{thm:MFPinq-geodesic} that 
\begin{align*}
  \Bigg(\fint_{B_d(x,r)}|u-u_{B_d(x,r)}|^q&d\mu\Bigg)^{1/q}\lesssim\left(\fint_{B_{\til d}(x,Lr)}|u-u_{B_{\til d}(x,Lr)}|^qd\mu\right)^{1/q}\\
  &\lesssim (1-\theta)r^\theta\fint_{B_{\til d}(x,Lr)}\int_{B_{\til d}(x,Lr)}\frac{|u(y)-u(z)|}{\til d(y,z)^\theta\mu(B_{\til d}(y,\til d(y,z)))}d\mu(z)d\mu(y)\\
  &\lesssim (1-\theta)r^\theta\fint_{B_{d}(x,L^2r)}\int_{B_{d}(x,L^2r)}\frac{|u(y)-u(z)|}{d(y,z)^\theta\mu(B_{d}(y,d(y,z)))}d\mu(z)d\mu(y).
\end{align*}
Thus, Theorem~\ref{thm: MFPinq} follows from Theorem~\ref{thm:MFPinq-geodesic}, with constants $C\ge 1$ and $\tau=L^2$ depending only on the original doubling and $(1,1)$-Poincar\'e inequality constants.

\begin{theorem}\label{thm:MFPinq-geodesic}
    Let $(X,d,\mu)$ be a complete, geodesic metric measure space, with $\mu$ a doubling measure supporting a $(1,1)$-Poincar\'e inequality.  Then there exists a constant $C\ge 1$ such that for all $0<\theta<1$, $1\le q\le Q_d/(Q_d-\theta)$, all balls $B_0:=B(x_0,r_0)\subset X$, and all $u\in L^1_\loc(X)$, we have 
    \[
    \left(\fint_{B_0}|u-u_{B_0}|^{q}d\mu\right)^{1/q}\le C(1-\theta)r_0^\theta\fint_{B_0}\int_{B_0}\frac{|u(x)-u(y)|}{d(x,y)^\theta\mu(B(x,d(x,y)))}d\mu(y)d\mu(x).
    \]
    Here, $Q_d>1$ is the relative lower mass bound exponent \eqref{eq:rel lower mass bound exponent}, and the constant $C$ depends only on the doubling and $(1,1)$-Poincar\'e inequality constants.
\end{theorem}

\begin{proof}
    Let $Q:=Q_d$, and let $q=Q_d/(Q_d-\theta)$.  It suffices to prove the result for this value of $q$, as the other cases then follow by H\"older's inequality.  Furthermore, without loss of generality, we may assume that $0<r_0<2\diam(X)$.  Indeed, in the case that $X$ is bounded and $r_0\ge2\diam(X)$, then $B_0=B(x_0,3\diam(X)/2)$.  Establishing the desired inequality for this ball with prescribed radius $3\diam(X)/2$ implies the inequality with radius $r_0$. 

For $\lambda>0$, let 
    \[
    \Omega_\lambda:=\{x\in B_0:|u-u_{B_0}|>\lambda\},
    \]
    and note that 
    \[
    \lambda^{q-1}\mu(\Omega_\lambda)^{(q-1)/q}\le\left(\int_0^\lambda\mu(\Omega_t)^{1/q}dt\right)^{q-1}.
    \]
    By Cavalieri's principle, we then have that 
    
    \begin{align}\label{eq:Frac PI Cavalieri}
        \bigg(\int_{B_0}|u-&u_{B_0}|^qd\mu\bigg)^{1/q}=\left(q\int_0^\infty\lambda^{q-1}\mu(B_0\cap\{|u-u_{B_0}|>\lambda\})d\lambda\right)^{1/q}\nonumber\\
        &=\left(q\int_0^\infty\lambda^{q-1}\mu(\Omega_\lambda)^{(q-1)/q}\mu(\Omega_\lambda)^{1/q}d\lambda\right)^{1/q}\nonumber\\
        &\le\left(q\int_0^\infty\left(\int_0^\lambda\mu(\Omega_t)^{1/q}dt\right)^{q-1}\mu(\Omega_\lambda)^{1/q}d\lambda\right)^{1/q}\nonumber\\
        &\le q^{1/q}\left(\int_0^\infty\mu(\Omega_t)^{1/q}dt\right)^{(q-1)/q}\left(\int_0^\infty\mu(\Omega_\lambda)^{1/q}d\lambda\right)^{1/q}\nonumber\\
        &\le q^{1/q}\int_0^\infty\mu(\Omega_\lambda)^{1/q}d\lambda\nonumber\\
        &\le q^{1/q}\left(\int_0^\infty\mu(B_0\cap\{u-u_{B_0}>\lambda\})^{1/q}\,d\lambda+\int_0^\infty\mu(B_0\cap\{u_{B_0}-u>\lambda\})^{1/q}\,d\lambda\right)\nonumber\\        
       &\le2\int_{u_{B_0}}^\infty\mu(B_0\cap\{u>\lambda\})^{1/q}\,d\lambda+2\int_{-\infty}^{u_{B_0}}\mu(B_0\cap\{u<\lambda\})^{1/q}\,d\lambda.
    \end{align}
    If $u$ is replaced with $-u$, then the two terms on the right-hand side of the above inequality switch roles, and so it suffices to estimate the first term.
    
    Denoting
    \begin{equation*}
        m_u:=m_u(B_0)=\inf\{t\in\R:\mu(B_0\cap\{u>t\})<\mu(B_0)/2\},
    \end{equation*}
    we have that 
    \begin{align}\label{eq:I+II}
        \int_{u_{B_0}}^\infty&\mu(B_0\cap\{u>\lambda\})^{1/q}\,d\lambda\nonumber\\
        &=\int_{u_{B_0}}^{\max\{m_u,u_{B_0}\}}\mu(B_0\cap\{u>\lambda\})^{1/q}\,d\lambda+\int_{\max\{m_u,u_{B_0}\}}^\infty\mu(B_0\cap\{u>\lambda\})^{1/q}\,d\lambda\nonumber\\
        &=:I+II.
    \end{align}

    To estimate $I$, we may assume that $u_{B_0}<m_u$.  By definition of $u_{B_0}$, it follows that 
    \[
    \int_{u_{B_0}}^\infty\mu(B_0\cap\{u>\lambda\})d\lambda=\int_{-\infty}^{u_{B_0}}\mu(B_0\cap\{u<\lambda\})d\lambda,
    \]
    and by the definition of $m_u$, we have that 
    \[
    \int_{u_{B_0}}^{m_u}\mu(B_0\cap\{u>\lambda\})d\lambda\ge \frac{\mu(B_0)}{2}(m_u-u_{B_0}).
    \]
    Hence, it follows that
    \begin{align}\label{eq:I-first}
        I\le\int_{u_{B_0}}^{m_u}\mu(B_0\cap\{u>\lambda\})^{1/q}\,d\lambda&\le(m_u-u_{B_0})\mu(B_0)^{1/q}\nonumber\\
        &\le\frac{2}{\mu(B_0)^{(q-1)/q}}\int_{u_{B_0}}^{m_u}\mu(B_0\cap\{u>\lambda\})d\lambda\nonumber\\
        &\le\frac{2}{\mu(B_0)^{(q-1)/q}}\int_{u_{B_0}}^{\infty}\mu(B_0\cap\{u>\lambda\})d\lambda\nonumber\\
        &=\frac{2}{\mu(B_0)^{(q-1)/q}}\int_{-\infty}^{u_{B_0}}\mu(B_0\cap\{u<\lambda\})d\lambda\nonumber\\
        &\le\frac{2}{\mu(B_0)^{(q-1)/q}}\int_{-\infty}^{m_u}\mu(B_0\cap\{u<\lambda\})d\lambda\nonumber\\
        &\le 2\int_{-\infty}^{m_u}\mu(B_0\cap\{u<\lambda\})^{1/q}\,d\lambda,
    \end{align}
    where we have again used the definition of $m_u$ to obtain the last inequality.

    For $\lambda<m_u$, set
    \begin{equation*}
        K_\lambda:=\left\lceil{\log_2\left(\left(\frac{\mu(B_0)}{\mu(B_0\cap\{u<\lambda\})}\right)^{1/Q}\right)}\right\rceil.
    \end{equation*}
    For $\lambda<m_u$ and $k\ge K_\lambda$, we then have that 
    \begin{align*}
        \frac{1}{2^{kQ}}\le\frac{\mu(B_0\cap\{u<\lambda\})}{\mu(B_0)}\le\frac{1}{2},
    \end{align*}
    and so applying Lemma~\ref{lem:Fractional Isoperimetric}, we obtain
    \begin{align*}
    \mu(B_0\cap\{u<\lambda\})^{\frac{Q-1}{Q}}\lesssim\frac{2^k}{\mu(B_0)^{1/Q}}\int_{B_0}\fint_{A_k(x)}|\chi_{\{u<\lambda\}}(y)-\chi_{\{u<\lambda\}}(x)|d\mu(y)d\mu(x),
    \end{align*}
    where 
    \[
    A_k(x):=B_0\cap B(x,2^{-k}r_0)\setminus B(x,2^{-k-1}r_0).
    \]
Multiplying $2^{-k(1-\theta)}$ to both sides and summing over $k\ge K_\lambda$, it then follows that 
\begin{align*}
    \sum_{k\ge K_\lambda}2^{-k(1-\theta)}\mu(B_0&\cap\{u<\lambda\})^{\frac{Q-1}{Q}}\\
    &\lesssim\sum_{k\ge K_\lambda}\frac{2^{k\theta}}{\mu(B_0)^{1/Q}}\int_{B_0}\fint_{A_k(x)}|\chi_{\{u<\lambda\}}(y)-\chi_{\{u<\lambda\}}(x)|d\mu(y)d\mu(x).
\end{align*}
By our choice of $K_\lambda$, we have that 
\begin{align*}
    \sum_{k\ge K_\lambda}2^{-k(1-\theta)}=\frac{2^{-K_\lambda(1-\theta)}}{1-2^{-(1-\theta)}}&\ge\frac{2^{-(1-\theta)}}{1-2^{-(1-\theta)}}\left(\frac{\mu(B_0\cap\{u<\lambda\})}{\mu(B_0)}\right)^{\frac{1-\theta}{Q}}\\
    &\ge\frac{1}{2(1-\theta)}\left(\frac{\mu(B_0\cap\{u<\lambda\})}{\mu(B_0)}\right)^{\frac{1-\theta}{Q}},
\end{align*}
    and so combining this with the previous expression, we obtain
    \begin{align*}
        \mu(B_0\cap\{u<\lambda\})^{1/q}\lesssim\frac{(1-\theta)}{\mu(B_0)^{(q-1)/q}}\sum_{k\ge K_\lambda}2^{k\theta}\int_{B_0}\fint_{A_k(x)}|\chi_{\{u<\lambda\}}(y)-\chi_{\{u<\lambda\}}(x)|d\mu(y)d\mu(x).
    \end{align*}

We claim that for all $x\in B_0$ and $y\in A_k(x)$, 
\begin{equation}\label{eq:annulus sub ball}
\mu(A_k(x))\simeq\mu(B(x,d(x,y))).
\end{equation}
Indeed, if $x\in B(x_0,r_0/2)$, then $B(x,2^{-k}r_0)\subset B_0$, and so, as $X$ is connected, there exists $z\in A_k(x)$ such that $d(x,z)=3(2^{-k}r_0)/4$. Setting $\wtil B:=B(z,2^{-k-2}r_0)$, we then have that $\wtil B\subset A_k(x)$.  If $x\in B_0\setminus B(x_0,r_0/2)$, then joining $x$ to $x_0$ by a geodesic $\gamma$, there exists $z$ in $\gamma$ such that $d(x,z)=3(2^{-k}r_0)/4$. Setting $\wtil B:=B(z,2^{-k-2}r_0)$, it follows that $\wtil B\subset A_k(x)$ in this case as well, since $\gamma$ is a geodesic.  For all $y\in A_k(x)$, it follows from doubling that in either case, we have
\begin{equation*}
\mu(\wtil B)\le\mu(A_k(x))\le\mu(B(x,2d(x,y)))\lesssim\mu(B(x,d(x,y)))\le\mu(8\wtil B)\lesssim\mu(\wtil B).
\end{equation*}
This proves the claim.

Applying this claim to the previous estimate, we have that 
\begin{align*}
    \mu(B_0\cap\{u<&\lambda\})^{1/q}\\
    &\lesssim\frac{(1-\theta)}{\mu(B_0)^{(q-1)/q}}\sum_{k\ge K_\lambda}2^{k\theta}(2^{-k}r_0)^\theta\int_{B_0}\int_{A_k(x)}\frac{|\chi_{\{u<\lambda\}}(y)-\chi_{\{u<\lambda\}}(x)|}{d(x,y)^\theta\mu(B(x,d(x,y)))}d\mu(y)d\mu(x)\\
        &=(1-\theta)\frac{r_0^\theta}{\mu(B_0)^{(q-1)/q}}\sum_{k\ge K_\lambda}\int_{B_0}\int_{A_k(x)}\frac{|\chi_{\{u<\lambda\}}(y)-\chi_{\{u<\lambda\}}(x)|}{d(x,y)^\theta\mu(B(x,d(x,y)))}d\mu(y)d\mu(x)\\
        &\le(1-\theta)\frac{r_0^\theta}{\mu(B_0)^{(q-1)/q}}\int_{B_0}\int_{B_0}\frac{|\chi_{\{u<\lambda\}}(y)-\chi_{\{u<\lambda\}}(x)|}{d(x,y)^\theta\mu(B(x,d(x,y)))}d\mu(y)d\mu(x).
\end{align*}
From \eqref{eq:I-first}, we then obtain 
\begin{align}\label{eq:I final}
    I&\lesssim\frac{(1-\theta)r_0^\theta}{\mu(B_0)^{(q-1)/q}}\int_{B_0}\int_{B_0}\left(\int_{-\infty}^{m_u}|\chi_{\{u<\lambda\}}(y)-\chi_{\{u<\lambda\}}(x)|d\lambda\right) \frac{d\mu(y)d\mu(x)}{d(x,y)^\theta\mu(B(x,d(x,y)))}\nonumber\\
    &\le(1-\theta) \frac{r_0^\theta}{\mu(B_0)^{-1/q}}\fint_{B_0}\int_{B_0}\frac{|u(x)-u(y)|}{d(x,y)^\theta\mu(B(x,d(x,y)))}d\mu(y)d\mu(x).
\end{align}

To estimate $II$, we first have that 
\begin{align}\label{eq:II-first}
    II\le\int_{m_u}^\infty\mu(B_0\cap\{u>\lambda\})^{1/q}\,d\lambda.
\end{align}
Using the same argument following \eqref{eq:I-first} used to estimate $I$, but replacing $\{u<\lambda\}$ for $\lambda<m_u$ with $\{u>\lambda\}$ for $\lambda>m_u$, we similarly obtain
\[
II\lesssim(1-\theta)\frac{r_0^\theta}{\mu(B_0)^{-1/q}}\fint_{B_0}\int_{B_0}\frac{|u(x)-u(y)|}{d(x,y)^\theta\mu(B(x,d(x,y)))}d\mu(y)d\mu(x).
\]
Combining this estimate with \eqref{eq:I final}, \eqref{eq:I+II}, and \eqref{eq:Frac PI Cavalieri} completes the proof.
\end{proof}

  As a corollary to Theorem~\ref{thm: MFPinq}, we obtain the following fractional relative isoperimetric inequality, which is a fractional analog of \eqref{eq:rel isoperimetric}.

\begin{proof}[Proof of Corollary \ref{thm: MFRiso}]
    The conclusion follows by setting $u:=\chi_E$ and applying Theorem~\ref{thm: MFPinq}:
    \begin{align*}
        \left(\frac{\min\{\mu(B\cap E),\,\mu(B\setminus E)\}}{\mu(B)}\right)^{1/q}&\le \left(2\frac{\mu(B\cap E)}{\mu(B)}\frac{\mu(B\setminus E)}{\mu(B)}\right)^{1/q}\\      &=\left(\fint_B\fint_B|u(x)-u(y)|^qd\mu(y)d\mu(x)\right)^{1/q}\\
        &\le \left(2^{q+1}\fint_B|u-u_B|^qd\mu\right)^{1/q}\\
        &\lesssim(1-\theta)r^\theta\fint_{\tau B}\int_{\tau B}\frac{|u(x)-u(y)|}{d(x,y)^\theta\mu(B(x,d(x,y)))}d\mu(y)d\mu(x)\\
        &\simeq(1-\theta)r^\theta\frac{P_\theta(E,\tau B)}{\mu(\tau B)}.\qedhere
    \end{align*}
\end{proof}

Finally, we give the other implication from the relative fractional perimeter inequality to the fractional Poincar\'e  inequality, which completes the proof of Theorem \ref{rmk: PiRiequiv}. 
    \begin{proof}[Proof of Theorem \ref{rmk: PiRiequiv}]
    The Corollary~\ref{thm: MFRiso} shows that if $(X,d,\mu)$ supports a $(\theta,q,1)_{BBM}$-Poincar\'e inequality, then it supports a $(\theta,q,1)_{BBM}$-relative isoperimetric inequality.  To prove the other direction, we fix a ball $B:=B(x,r)\subset X$ and $u\in L^1_\loc(X)$. As in the proof of Theorem \ref{thm: MFPinq}, it suffices to estimate the sum 
    \[
        I+II:=\int_{u_{B}}^{\max\{m_u,u_{B}\}}\mu(B\cap\{u>\lambda\})^{1/q}\,d\lambda+\int_{\max\{m_u,u_{B}\}}^\infty\mu(B\cap\{u>\lambda\})^{1/q}\,d\lambda.
    \]
     Since $\mu(B\cap\{u<\lambda\})\leq \mu(B)/2$ for $\lambda<m_u$, it follows from \eqref{eq:I-first}, inequality \eqref{MFRiso}, and the fractional coarea formula Lemma~\ref{lem:frac coarea} that 
    \begin{align}\label{eq:estimate-I}
    I&\le2\int_{-\infty}^{m_u}\mu(B\cap\{u<\lambda\})^{1/q}d\lambda\nonumber\\
    &\lesssim (1-\theta)\frac{r^\theta}{\mu(B)^{(q-1)/q}}\int_{-\infty}^{m_u} P_\theta(\{u<\lambda\},\tau B)d\lambda\nonumber\\
    &\lesssim (1-\theta)\frac{r^\theta}{\mu(B)^{(q-1)/q}}\int_{\R} P_\theta(\{-u>\lambda\},\tau B)d\lambda\nonumber\\
    &\lesssim(1-\theta)\frac{r^\theta}{\mu(B)^{-1/q}}\fint_{\tau B}\int_{\tau B}\frac{|u(x)-u(y)|}{d(x,y)^\theta\mu(B(x,d(x,y)))}d\mu(y)d\mu(x).
    \end{align}
    
    To estimate $II$, we note that $\mu(B\cap\{u>\lambda\})\leq \mu(B)/2$ for $\lambda>m_u$.  Similar to the estimate for $I$, we then obtain 
    \begin{align}\label{eq:estimate-II}
        II\lesssim(1-\theta)\frac{r^\theta}{\mu(B)^{-1/q}}\fint_{\tau B}\int_{\tau B}\frac{|u(x)-u(y)|}{d(x,y)^\theta\mu(B(x,d(x,y)))}d\mu(y)d\mu(x).
    \end{align}
    Combining \eqref{eq:estimate-I} and \eqref{eq:estimate-II} with \eqref{eq:Frac PI Cavalieri} and \eqref{eq:I+II} then yields 
 \[
    \left(\fint_{B}|u-u_{B}|^{q}d\mu\right)^{1/q}\lesssim (1-\theta)r^\theta\fint_{\tau B}\int_{\tau B}\frac{|u(x)-u(y)|}{d(x,y)^\theta\mu(B(x,d(x,y)))}d\mu(y)d\mu(x).\qedhere
    \]    
    \end{proof}

\section{Equivalence between Poincar\'e inequalities}\label{FPIeqv}

In this section, we prove Theorem~\ref{fractional to Poincare}, which shows that if $(X,d,\mu)$ supports a $(\theta,q,1)_{BBM}$-Poincar\'e inequality for $\theta\to 1^-$, then it supports a $(q,1)$-Poincar\'e inequality.  First we essentially repeat some arguments from \cite{LPZ}. 

Let $\Om\subset X$ denote an open set.
Let $u$ be a function defined on $\Om$. We define
\begin{equation}\label{eq:Lip r}
	\Lip_r u(x):=\sup_{y\in \Om\cap B(x,r) \setminus \{x\}}\frac{|u(y)-u(x)|}{d(y,x)},
	\quad x\in \Om,\quad r>0,
\end{equation}
and
\begin{equation}\label{eq:Lip}
	\Lip u(x):=\limsup_{r\to 0}\Lip_ru(x).
\end{equation}

We consider the following kernels.
We consider a sequence of nonnegative {$\mu\times \mu$}-measurable functions
$\{\rho_i(x,y)\}_{i=1}^{\infty}$, $x,y\in X$,
and a fixed constant
$1\le C_{\rho}<\infty$ satisfying the following conditions:
\begin{enumerate}[(1)]\label{rho conditions}
	\item For every $x,y\in X$ with $0<d(x,y)< 1$, we have
	\begin{equation}\label{eq:rho hat majorize}
		\sum_{j=1}^{\infty}d_{i,j}'\frac{\chi_{B(y,2^{-j+1})
				\setminus B(y,2^{-j})}(x)}{\mu(B(y,2^{-j+1}))}
		\le \rho_i(x,y)
		\le \sum_{j=1}^{\infty}d_{i,j}\frac{\chi_{B(y,2^{-j+1})
				\setminus B(y,2^{-j})}(x)}{\mu(B(y,2^{-j+1}))}
	\end{equation}
	for numbers $d_{i,j}',d_{i,j}\ge 0$ for which 
	$\sum_{j=1}^{\infty} d_{i,j}'\ge C_\rho^{-1}$ and
	$\sum_{j=1}^{\infty} d_{i,j}\le C_\rho$ for all $i\in\N$.
	
	\item For all $\delta>0$, we have
	\begin{equation}\label{eq:radius one condition}
		\lim_{i\to\infty}	\left(\sup_{y\in \Om}\int_{\Om\setminus B(y,\delta)}\frac{\rho_i(x,y)}{d(x,y)}\,d\mu(x)
		+\sup_{x\in \Om}\int_{\Om\setminus B(x,\delta)}\frac{\rho_i(x,y)}{d(x,y)}\,d\mu(y)\right)=0.
	\end{equation}
\end{enumerate}

In particular, we can take
\[
\rho_i(x,y)=\frac{1-\theta_i}{d(x,y)^{\theta_i-1}\mu(B(x,d(x,y)))},
\]
for which we can choose $d'_{i,j}=(1-\theta_i)2^{-j(1-\theta_i)}$ and
$d_{i,j}=C_\mu (1-\theta_i)2^{(-j+1)(1-\theta_i)}$ {where $\theta_i\in (0,1)$ is a sequence converging to $1$}.
Indeed, for every $x,y\in X$ with $0<d(x,y)<1$,  we now have
\begin{align*}
	\sum_{j=1}^{\infty}d'_{i,j} \frac{\chi_{B(y,2^{-j+1})
			\setminus B(y,2^{-j})}(x)}{\mu(B(y,2^{-j+1}))} 
	&\le \frac{1-\theta_i}{d(x,y)^{\theta_i-1}\mu(B(y,d(x,y)))}\\
	&\le \sum_{j=1}^{\infty}d_{i,j} \frac{\chi_{B(y,2^{-j+1})
			\setminus B(y,2^{-j})}(x)}{\mu(B(y,2^{-j+1}))}.
\end{align*}
Here
\[
\sum_{j=1}^{\infty}2^{-j(1-\theta_i)}
=\frac{2^{\theta_i-1}}{1-2^{\theta_i-1}}
\ge \frac{2^{\theta_i-1}}{1-\theta_i}
\ge \frac{1/2}{1-\theta_i},
\]
and so
\[
\sum_{j=1}^{\infty}d'_{i,j}\ge 1/2,
\]
satisfying the first inequality given after \eqref{eq:rho hat majorize}.
Similarly,
\[
\sum_{j=1}^{\infty}2^{(-j+1)(1-\theta_i)}\le 2\int_0^2 t^{(1-\theta_i)-1} dt=\frac{2^{2-\theta_i}}{1-\theta_i}\le \frac{4}{1-\theta_i},
\]
and so
\[
\sum_{j=1}^{\infty}d_{i,j}\le 4C_\mu,
\]
satisfying  the second inequality given after \eqref{eq:rho hat majorize}.
Finally, we estimate

\begin{align*}
	\frac{\rho_i(x,y)}{d(x,y)}
	&=(1-\theta_i)\frac{1}{d(x,y)^{ \theta_i}\mu(B(x,d(x,y)))}\\
	&\le C_\mu(1-\theta_i)\sum_{j\in \Z}2^{j\theta_i} \frac{\chi_{B(y,2^{-j+1})
			\setminus B(y,2^{-j})}(x)}{\mu(B(y,2^{-j+1}))},
\end{align*}
and so (the notation $\sum_{j\le  -\log_2 \delta}$ means that we sum over
integers $j$ at most $\lceil-\log_2 \delta\rceil$)
\begin{align*}
	\int_{X\setminus B(y,\delta)}\frac{\rho_i(x,y)}{d(x,y)}\,d\mu(x)
	&\le C_\mu(1-\theta_i) \sum_{j\le  -\log_2 \delta}2^{j\theta_i}\\
	&\le C_\mu (1-\theta_i)\frac{\delta^{-\theta_i}}{1-2^{-\theta_i}}\\
	&\to 0\quad\textrm{as }{\theta_i\to 1 \text{ when  }i\to \infty},
\end{align*}
and so \eqref{eq:radius one condition} holds.

\begin{lemma}\label{prop:one direction}
	Let $u\in \Lip(X)$ and $0<R\le 1$, and suppose $U\subset X$ is open.
	Suppose $\{\rho_i\}_{i=1}^{\infty}$ is a sequence of mollifiers that satisfy
	\eqref{eq:rho hat majorize}.
	Then
	\begin{equation}\label{eq:upper bound with 8 lambda}
		\int_{U} \int_{B(y,R)} \frac{|u(x)-u(y)|}{d(x,y)}\rho_i(x,y)\,d\mu(x)\,d\mu(y)
		\le C_{\rho}\int_{U} \Lip_{R}u(y) \,d\mu(y)
	\end{equation}
	for every $i\in\N$.
	\end{lemma}

\begin{proof}
	We have
	\begin{align*}
		\int_{U} \int_{B(y,R)} \frac{|u(x)-u(y)|}{d(x,y)}\rho_i(x,y)\,d\mu(x)\,d\mu(y)
		&\le \int_{U} \int_{B(y,R)} \Lip_{R}u(y) \rho_i(x,y)\,d\mu(x)\,d\mu(y)\\
		&\le C_{\rho}\int_{U} \Lip_{R}u(y)\,d\mu(y)
	\end{align*}
	by \eqref{eq:rho hat majorize}.
\end{proof}

\begin{theorem}\label{thm:BBM one direction}
	Suppose $\Om\subset X$ is open and bounded, 
	and let $u\in \Lip(\Om)$.
	Suppose $\{\rho_i\}_{i=1}^{\infty}$ is a sequence of mollifiers that satisfy
	\eqref{eq:rho hat majorize} and \eqref{eq:radius one condition}.
	Then
	\[
	\limsup_{i\to\infty}\int_{\Om}\int_{\Om} \frac{|u(x)-u(y)|}{d(x,y)}\rho_i(x,y)\,d\mu(x)\,d\mu(y)
	\le C_{\rho}\int_{\Om} \Lip u \,d\mu.
	\]
\end{theorem}

\begin{proof}
	Consider $0<R\le 1$. By a standard extension, we can assume that $u\in \Lip(X)$. We have
	\begin{align*}
		&\int_{\Om} \int_{\Om}\frac{|u(x)-u(y)|}{d(x,y)}\rho_i(x,y)\,d\mu(x)\,d\mu(y)\\
		&\qquad\le  \int_{\Om} \int_{\Om\setminus B(y,R)}\frac{|u(x)-u(y)|}{d(x,y)}\rho_i(x,y)\,d\mu(x)\,d\mu(y)\\
		&\qquad\qquad + \int_{\Om} \int_{ B(y,R)}\frac{|u(x)-u(y)|}{d(x,y)}\rho_i(x,y)\,d\mu(x)\,d\mu(y).
	\end{align*}
	For the first term, we estimate 
	\begin{align*}
		&\int_{\Om} \int_{\Om\setminus B(y,R)}\frac{|u(x)-u(y)|}{d(x,y)}\rho_i(x,y)\,d\mu(x)\,d\mu(y)\\
		&\qquad \le 2\int_{\Om} \int_{\Om\setminus B(y,R)}\frac{|u(x)|+|u(y)|}{d(x,y)}\rho_i(x,y)\,d\mu(x)\,d\mu(y)\\
		&\qquad \le 2\int_{\Om}|u(y)| \int_{\Om\setminus B(y,R)}\frac{\rho_i(x,y)}{d(x,y)}\,d\mu(x)\,d\mu(y)\\
		&\qquad \qquad +2\int_{\Om}|u(x)| \int_{\Om\setminus B(x,R)}\frac{\rho_i(x,y)}{d(x,y)}\,d\mu(y)\,d\mu(x)\quad\textrm{by Fubini}\\
		&\qquad \le 2\int_{\Om}|u|\,d\mu 
		\left(\sup_{y\in \Om}\int_{\Om\setminus B(y,R)}\frac{\rho_i(x,y)}{d(x,y)}\,d\mu(x)
		+\sup_{x\in \Om}\int_{\Om\setminus B(x,R)}\frac{\rho_i(x,y)}{d(x,y)}\,d\mu(y)\right)\\
		&\qquad \to 0
	\end{align*}
	as $i\to\infty$ by \eqref{eq:radius one condition}.
	
	For the second term, from Lemma \ref{prop:one direction} we get for every $i\in\N$
	\begin{align*}
		\int_{\Om} \int_{B(y,R)}\frac{|u(x)-u(y)|}{d(x,y)}\rho_i(x,y)\,d\mu(x)\,d\mu(y)
		&\le C_{\rho} \int_{\Om} \Lip_{R} u \,d\mu\\
	\end{align*}

	Passing first to the $\limsup_{i\to\infty}$ in the two terms and then sending $R\to 0$, we get,
	using also Lebesgue's dominated convergence theorem
	\[
	\limsup_{i\to\infty}\int_{\Om} \int_{\Om}\frac{|u(x)-u(y)|}{d(x,y)}\rho_i(x,y)\,d\mu(x)\,d\mu(y)
	\le C_{\rho}\int_{\Om} \Lip u  \,d\mu.
	\]
\end{proof}

\begin{theorem}\label{thm:rho and Poincare}
	Suppose $\{\rho_i\}_{i=1}^{\infty}$ is a sequence of mollifiers that satisfy
	\eqref{eq:rho hat majorize} and \eqref{eq:radius one condition}, and suppose
	$\theta_i\in (0,1)$ with $\theta_i\to 1$
	as $i\to\infty$, and let $1\le q<\infty$.
	Suppose there exist constants $C,\tau \ge 1$ such that for all
	$i\in\N$,  all balls $B:=B(x,r)\subset X$, and all $u\in L^1_\loc(X)$, we have 
	\[
	\left(\fint_{B}|u-u_{B}|^q\,d\mu\right)^{1/q}
	\le Cr^{\theta_i}\fint_{\tau B}\int_{\tau B}\frac{|u(x)-u(y)|}{d(x,y)}\rho_i(x,y)\,d\mu(x)\,d\mu(y).
	\]
	Then $X$ supports a $(q,1)$-Poincar\'e inequality.
\end{theorem}
\begin{proof}
	Let $u\in \Lip(X)$. Consider a ball $B(x,r)$.  By Theorem \ref{thm:BBM one direction} with the choice $\Om=B(x,\tau r)$ we have 
	\[
	\limsup_{i\to\infty}\fint_{B(x,\tau r)}\int_{B(x,\tau  r)}
	\frac{|u(x)-u(y)|}{d(x,y)}\rho_i(x,y)\,d\mu(x)\,d\mu(y)
	\le C_{\rho}{\fint_{B(x,\tau r)}} \Lip u \,d\mu.
	\]
	It follows that
	\[
	\left(\fint_{B(x,r)}|u-u_{B(x,r)}|^q\,d\mu\right)^{1/q}
	\le C C_{\rho} r{\fint_{B(x,\tau r)}} \Lip u \,d\mu.
	\]
 As this inequality holds for all $u\in\Lip(X)$, it follows from \cite[Theorem~3.1]{KoLa14} that $X$ supports a $(q,1)$-Poincar\'e inequality.
\end{proof}

Theorem~\ref{fractional to Poincare} is then proven with the following corollary:

\begin{corollary}\label{cor:fractional to Poincare}
	Let $1\le q<\infty$ and suppose there exist constants $C,\tau \ge 1$ and a sequence of numbers $\theta_i\in (0,1)$ converging to $1$, such that for all balls $B:=B(x,r)\subset X$, and all $u\in L^1_\loc(X)$, we have 
	\[
	\left(\fint_{B}|u-u_{B}|^q\,d\mu\right)^{1/q}
	\le C(1-\theta_i)r^{\theta_i}\fint_{\tau B}\int_{\tau B}\frac{|u(x)-u(y)|}{d(x,y)^{\theta_i}\mu(B(x,d(x,y)))}\,d\mu(x)\,d\mu(y).
	\]
	Then $X$ supports a $(q,1)$-Poincar\'e inequality.
\end{corollary}
\begin{proof}
	Apply Theorem \ref{thm:rho and Poincare} with the kernels
	\[
	\rho_i(x,y)=\frac{1-\theta_i}{d(x,y)^{\theta_i-1}\mu(B(x,d(x,y)))}
	\]
	with $\theta_i\to 1$ as $i\to\infty$. The conditions
	\eqref{eq:rho hat majorize} and \eqref{eq:radius one condition} were verified for these kernels
	after they were introduced.
\end{proof}

\section{Fractional boxing inequality}\label{FBoxing}

In this section we prove the fractional boxing inequality Theorem~\ref{MFBoxing}, as well as a local variant Proposition~\ref{thm:loc-box-ineq}, under the assumptions of Theorem~\ref{thm: MFPinq}.  Recall that as $(X,d,\mu)$ is doubling and connected under these assumptions, the reverse doubling condition holds. That is, there exist $s>0$ and $C_s\ge 1$, depending only on the doubling constant, such that 
\begin{equation}\label{eq:reverse doubling}
\frac{\mu(B(y,r))}{\mu(B(x,R))}\le C_s\left(\frac{r}{R}\right)^s
\end{equation}
for all $0<r\le R\le 2\diam X$, $x\in X$ and $y\in B(x,R)$.

We begin with the following lemma, which concerns the large-scale behavior of the fractional perimeter functional.  This is a generalization to the metric setting of \cite[Lemma~4.2]{PS20}.
\begin{lemma}\label{lem:theta frac iso}
    Let $0<\theta<1$, and let $(X,d,\mu)$ be doubling and connected with $\mu(X)=\infty$. Given $0<\gamma<1$, there exists a constant $C\ge 1$, depending only on $\gamma$ and the doubling constant such that 
    \begin{align*}
        \frac{\mu(B(x_0,R))}{R^{\theta}}\le C\theta\int_{B(x_0,R)\cap E}\int_{X\setminus E}\frac{1}{d(x,y)^{\theta}\mu(B(x,d(x,y)))}d\mu(y)d\mu(x)
    \end{align*}
    for all measurable sets $E\subset X$, $x_0\in X$, and $R>0$ satisfying
    \begin{align}\label{eq:growth density}
        \frac{\mu(B(x_0,2^{-1}R)\cap E)}{\mu(B(x_0,2^{-1}R))}\ge \gamma\quad\text{ and }\quad\frac{\mu(B(x_0,r)\cap E)}{\mu(B(x_0,r))}<\gamma~\text{ for all }r\ge R.
    \end{align}
\end{lemma}

\begin{proof}
    If $x\in B(x_0,R)$ and $y\in X\setminus B(x_0,2R)$, then $d(x,y)\le 2d(x_0,y)$.  Thus, by \eqref{eq:growth density} and doubling, we have that 
    \begin{align*}\label{eq:R^Q_1}
        \int_{B(x_0,R)\cap E}&\int_{X\setminus E}\frac{1}{d(x,y)^{\theta}\mu(B(x,d(x,y)))}d\mu(y)d\mu(x)\\\nonumber
        &\ge\int_{B(x_0,R)\cap E}\int_{X\setminus(B(x_0,2R)\cup E)}\frac{1}{d(x,y)^{\theta}\mu(B(x,d(x,y)))}d\mu(y)d\mu(x)\nonumber\\
        &\gtrsim\mu(B(x_0,R)\cap E)\int_{X\setminus(B(x_0,2R)\cup E)}\frac{1}{d(x_0,y)^{\theta}\mu(B(x_0,d(x_0,y)))}d\mu(y)\nonumber\\
        &\ge\mu(B(x_0,R/2)\cap E)\int_{X\setminus(B(x_0,2R)\cup E)}\frac{1}{d(x_0,y)^{\theta}\mu(B(x_0,d(x_0,y)))}d\mu(y)\nonumber\\
        &\ge\gamma\mu(B(x_0,R/2))\int_{X\setminus(B(x_0,2R)\cup E)}\frac{1}{d(x_0,y)^{\theta}\mu(B(x_0,d(x_0,y)))}d\mu(y)\nonumber\\
        &\gtrsim \gamma\mu(B(x_0,R))\int_{X\setminus(B(x_0,2R)\cup E)}\frac{1}{d(x_0,y)^{\theta}\mu(B(x_0,d(x_0,y)))}d\mu(y).\nonumber
    \end{align*}
    Setting 
    \[
    \lambda_0:=\frac{1}{2}\left(\frac{1-\gamma}{2C_s}\right)^{1/s}, 
    \]
   where $s,C_s$ come from the reverse doubling condition  \eqref{eq:reverse doubling} and  
    \[
    A_i:=B(x_0,2R\lambda_0^{-i-1})\setminus B(x_0,2R\lambda_0^{-i})
    \]
    for $i\in\N\cup\{0\}$, we then have
     \begin{align}\label{eq:R^Q 1}
        \int_{B(x_0,R)\cap E}\int_{X\setminus E}&\frac{1}{d(x,y)^{\theta}\mu(B(x,d(x,y)))}d\mu(y)d\mu(x)\nonumber\\
        &\gtrsim \gamma\mu(B(x_0,R))\sum_{i=0}^\infty\int_{A_i\setminus E}\frac{1}{d(x_0,y)^\theta\mu(B(x_0,d(x_0,y)))}d\mu(y)\nonumber\\
        &\gtrsim\gamma R^Q\mu(B(x_0,R))\sum_{i=0}^\infty\frac{\lambda_0^{-iQ}}{\mu(B(x_0,2R\lambda_0^{-i-1}))}\int_{A_i\setminus E}\frac{1}{d(x_0,y)^{Q+\theta}}d\mu(y).
    \end{align}
    Here, the last inequality follows from the fact that $d(x_0,y)\simeq R\lambda_0^{-i}$ for all $y\in A_i$.

    Setting 
    \[    f_i(y):=\frac{\chi_{A_i\setminus E}(y)}{d(x_0,y)},
    \]
    it follows that for $\lambda>0$, 
    \[
    \{y\in X:f_i(y)>\lambda\}=A_i\cap B(x_0,1/\lambda)\setminus E.
    \]
    Applying Cavalieri's principle, we obtain
    \begin{align*}
    \int_{A_i\setminus E}\frac{1}{d(x_0,y)^{Q+\theta}}d\mu(y)&=(Q+\theta)\int_0^{\lambda_0^i/(2R)}\lambda^{Q+\theta-1}\mu(A_i\cap B(x_0,1/\lambda)\setminus E)d\lambda\\
    &\ge Q\int^{\lambda_0^{i+1}/R}_{\lambda_0^{i+1}/(2R)}\lambda^{Q+\theta-1}\mu(B(x_0,1/\lambda)\setminus(B(x_0,2R\lambda_0^{-i})\cup E))d\lambda.
    \end{align*}
    For $\lambda_0^{i+1}/(2R)\le\lambda\le\lambda_0^{i+1}/R$, it follows from \eqref{eq:growth density}, \eqref{eq:reverse doubling}, and our choice of $\lambda_0$ that 
    \begin{align*}
        \mu(B(x_0,1/\lambda)\setminus(B(x_0,2R\lambda_0^{-i})\cup E))&\ge\mu(B(x_0,1/\lambda)\setminus B(x_0,2R\lambda_0^{-i}))-\mu(B(x_0,1/\lambda)\cap E)\\
        &\ge\mu(B(x_0,1/\lambda))\left(1-\frac{\mu(B(x_0,2R\lambda_0^{-i}))}{\mu(B(x_0,1/\lambda))}\right)-\gamma\mu(B(x_0,1/\lambda))\\
        &\ge(1-C_s(2R\lambda/\lambda_0^{i})^s-\gamma)\mu(B(x_0,1/\lambda))\\
        &\ge(1-C_s(2\lambda_0)^s-\gamma)\mu(B(x_0,1/\lambda))\\
        &=\frac{1-\gamma}{2}\mu(B(x_0,1/\lambda)).
    \end{align*}
    Therefore, by the previous expression, \eqref{eq:R^Q 1}, and \eqref{eq:rel lower mass bound exponent}, we have that 
    \begin{align*}
        \int_{B(x_0,R)\cap E}\int_{X\setminus E}&\frac{1}{d(x,y)^{\theta}\mu(B(x,d(x,y)))}d\mu(y)d\mu(x)\\
        &\gtrsim\gamma(1-\gamma) R^Q\mu(B(x_0,R))\sum_{i=0}^\infty\lambda_0^{-iQ}\int_{\lambda_0^{i+1}/(2R)}^{\lambda_0^{i+1}/(R)}\lambda^{Q+\theta-1}\frac{\mu(B(x_0,1/\lambda))}{\mu(B(x_0,2R\lambda_0^{-i-1}))}d\lambda\\
        &\gtrsim\gamma(1-\gamma) R^Q\mu(B(x_0,R))\sum_{i=0}^\infty\lambda_0^{-iQ}\int_{\lambda_0^{i+1}/(2R)}^{\lambda_0^{i+1}/(R)}\lambda^{Q+\theta-1}\left(\frac{\lambda_0^{i+1}}{2R\lambda}\right)^{Q}d\lambda\\
        &\gtrsim\gamma(1-\gamma)\mu(B(x_0,R))\sum_{i=0}^\infty\int_{\lambda_0^{i+1}/(2R)}^{\lambda_0^{i+1}/(R)}\lambda^{\theta-1}d\lambda\\
        &=\gamma(1-\gamma)\lambda_0^\theta(1-2^{-\theta})\frac{\mu(B(x_0,R))}{\theta R^\theta}\sum_{i=0}^\infty\lambda_0^{i\theta}\\
        &=\gamma(1-\gamma)\lambda_0^\theta\frac{(1-2^{-\theta})}{(1-\lambda_0^\theta)}\frac{\mu(B(x_0,R))}{\theta R^\theta}.
    \end{align*}
This yields the desired result, since 
\[
\lim_{\theta\to 0}\frac{1-2^{-\theta}}{1-\lambda_0^\theta}=\frac{\log2}{\log(1/\lambda_0)}.\qedhere
\]
\end{proof}

We now prove the main result of this section, Theorem~\ref{MFBoxing}:

\begin{proof}[Proof of Theorem~\ref{MFBoxing}]
    Fix $x_0\in U$.  As $U$ is bounded and $\mu(X)=\infty$, there exists $R_0>0$ such that $U\subset B(x_0,R_0)$ and 
    \begin{equation*}
        \frac{\mu(B(x_0,R_0)\cap U)}{\mu(B(x_0,R_0))}\le 1/2.
    \end{equation*}
    For each $x\in U$, we have that $B(x_0,R_0)\subset B(x,2R_0)$. Therefore, for all $r\ge 2R_0$, we have that 
    \[
    \frac{\mu(B(x,r)\cap U)}{\mu(B(x,r))}\le\frac{\mu(B(x_0,R_0)\cap U)}{\mu(B(x_0,R_0))}\le1/2.
    \]
    Moreover, by the Lebesgue differentiation theorem, it follows that there exists $\mathcal{N}\subset U$ with $\mu(\mathcal{N})=0$ such that for all $x\in U\setminus\mathcal{N}$, we have that
    \[
	\lim_{r\to 0^+}\frac{\mu(B(x,r)\cap U)}{\mu(B(x,r))}=1.
   	\]
	Setting 
    \[
    \alpha_x:=\sup\{r>0:\mu(B(x,r)\cap U)\ge \mu(B(x,r))/2\}
    \]
	for $x\in U\setminus\mathcal{N}$,
    it then follows that $0<\alpha_x\le 2R_0$. We then choose $R_x>\alpha_x$ such that $R_x/2\le\alpha_x$ and 
    \begin{align}\label{eq:>1/2}
        \frac{\mu(B(x,R_x/2)\cap U)}{\mu(B(x,R_x/2))}\ge 1/2.
    \end{align}
    Since $R_x>\alpha_x$, we also have that 
\begin{align}\label{eq:<1/2}
        \frac{\mu(B(x,r)\cap U)}{\mu(B(x,r))}<1/2
    \end{align}
    for all $r\ge R_x$.  Note that $R_x\le 4R_0$. 
    
    By doubling, \eqref{eq:>1/2}, and \eqref{eq:<1/2}, it follows that
    \begin{align}\label{eq:Rx density}
        \frac{1}{2}>\frac{\mu(B(x,R_x)\cap U)}{\mu(B(x,R_x))}\ge\frac{\mu(B(x,2^{-1}R_x)\cap U)}{C_\mu\mu(B(x,2^{-1}R_x))}\ge\frac{1}{2C_\mu}.
    \end{align}
    From this inequality and applying the $(\theta,1,1)_{BBM}$-relative isoperimetric inequality given by Corollary~\ref{thm: MFRiso}, we obtain
    \begin{align*}
        \frac{1}{2C_\mu}\mu(B(x,R_x))\le\mu(B(x,R_x)\cap U)&=\min\{\mu(B(x,R_x)\cap U),\,\mu(B(x,R_x)\setminus U)\}\\
        &\le C(1-\theta)R_x^\theta P_\theta(U, B(x,\tau R_x)),
    \end{align*}
    and so it follows that 
    \begin{align}\label{eq:Rx}
        \frac{\mu(B(x,R_x))}{R_x^\theta}\lesssim (1-\theta)\int_{ B(x,\tau R_x)\cap U}\int_{B(x,\tau R_x)\setminus U}\frac{1}{d(z,y)^\theta\mu(B(z,d(z,y)))}d\mu(y)d\mu(z).
    \end{align}
    
    As the collection $\{B(x,\tau R_x)\}_{x\in U\setminus\mathcal{N}}$ has uniformly bounded radii, we apply the 5-covering lemma to obtain a countable, pairwise disjoint subcollection $\{B(x_i,\tau r_i)\}_{i\in\N}$ such that 
    \begin{align*}
        U\setminus\mathcal{N}\subset\bigcup_{i\in\N}B(x_i,5\tau r_i).
    \end{align*}
    By pairwise disjointness of this subcollection and \eqref{eq:Rx}, we then obtain
    \begin{align}\label{eq:1-theta boxing}
        \sum_{i\in\N}\frac{\mu(B(x_i,r_i))}{r_i^\theta}&\lesssim (1-\theta)\sum_{i\in\N}\int_{B(x_i,\tau r_i)\cap U}\int_{B(x_i,\tau r_i)\setminus U}\frac{1}{d(x,y)^\theta\mu(B(x,d(x,y)))}d\mu(y)d\mu(x)\nonumber
        \\
        &\le(1-\theta)\int_{U}\int_{X\setminus U}\frac{1}{d(x,y)^\theta\mu(B(x,d(x,y)))}d\mu(y)d\mu(x).
    \end{align}
   
    The estimates \eqref{eq:>1/2} and \eqref{eq:<1/2} also allow us to apply Lemma~\ref{lem:theta frac iso}, with $E=U$, $x_0=x_i$, $R=r_i$, and $\gamma=1/2$, to obtain  
    \begin{align*}
        \frac{\mu(B(x_i,r_i))}{r_i^\theta}\lesssim\theta\int_{B(x_i,r_i)\cap U}\int_{X\setminus U}\frac{1}{d(x,y)^\theta\mu(B(x,d(x,y)))}d\mu(y)d\mu(x).
    \end{align*}
    Summing over $i\in\N$, we then have
    \begin{align*}
        \sum_{i\in\N} \frac{\mu(B(x_i,r_i))}{r_i^\theta}\lesssim \theta P_\theta(U,X).
    \end{align*}
    Since $\min\{\theta,1-\theta\}\le 2\theta(1-\theta)$, combining this inequality with \eqref{eq:1-theta boxing} yields \eqref{eq:doubling boxing ineq}.    
\end{proof}

\begin{remark}\label{rem:open set boxing}
We note that if $U$ is additionally assumed to be open, then 
\[
\lim_{r\to 0^+}\frac{\mu(B(x,r)\cap U)}{\mu(B(x,r))}=1
\]
holds for all $x\in U$, and so the conclusion of Theorem~\ref{MFBoxing} holds without requiring the negligible set $\mathcal{N}$. 
\end{remark}

\begin{remark}\label{Remark-lower-boxing}
If we assume in addition that there exists $t>0$ and $C_t\ge 1$ such that 
    \begin{equation}\label{eq:lower mass bound}
    C_t^{-1}R^t\le\mu(B(x,R))
    \end{equation}
    for all $x\in X$ and $0<R<2\diam(X)$, i.e.\ $\mu$ is lower Ahlfors $t$-regular, then the conclusion of Theorem~\ref{MFBoxing} holds under the weaker hypothesis that $\mu(U)<\infty$.  Indeed, under this assumption, we have by \eqref{eq:Rx density} in the above proof that 
    \begin{align*}
        R_x^t\le C_t\mu(B(x,R_x))\le 2C_tC_\mu\mu(B(x,R_x)\cap U)\le 2C_tC_\mu\mu(U)<\infty.
    \end{align*}
    Thus, the radii $\{R_x\}_{x\in U}$ are uniformly bounded even if $U$ is unbounded.  This allows us to apply the $5$-covering lemma, and so the proof proceeds in the same manner.

   If instead of \eqref{eq:lower mass bound}, we additionally assume that there exists $c_0>0$ such that 
\begin{equation}\label{eq:unit ball lower mass bound}
\mu(B(x,1))\ge c_0
\end{equation}
for all $x\in X$, then the radii $R_x$ are similarly uniformly bounded, and so Theorem~\ref{MFBoxing} also holds under the weaker hypothesis that $\mu(U)<\infty$.  To see this, we note that if $x\in X$ and $0<r<1$, then by the doubling assumption and \eqref{eq:unit ball lower mass bound}, we have that 
\[
\mu(B(x,r))\ge C_Q^{-1}r^{Q_d}\mu(B(x,1))\ge C_Q^{-1}c_0 r^{Q_d},
\]
where $Q_d$ is as in \eqref{eq:rel lower mass bound exponent} and $C_Q$ depends only on the doubling constant.  If $x\in X$ and $r>1$, then by \eqref{eq:reverse doubling}, which is guaranteed by the assumptions of doubling and the $1$-Poincar\'e inequality, and \eqref{eq:unit ball lower mass bound}, we have that 
\[
\mu(B(x,r))\ge C_s^{-1}r^s\mu(B(x,1))\ge C_s^{-1}c_0 r^s.
\]
Hence,  for all $x\in X$ and $r>0$, it follows that
 \[
\mu(B(x,r))\ge c_0\min\{C_s^{-1}r^s, C_Q^{-1}r^{Q_d}\}.
\]
Considering the radii $\{R_x\}_{x\in U}$ from the proof of Theorem~\ref{MFBoxing}, it follows from \eqref{eq:Rx density} that $\mu(B(x,R_x))\le 2C_\mu\mu(U)$.  By the previous inequality, it then follows that 
\[
R_x\le c_0^{-1}2C_\mu\max\{(C_s\mu(U))^{1/s},(C_Q\mu(U))^{1/{Q_d}}\}<\infty.
\]


    If $(X,d,\mu)$ satisfies neither \eqref{eq:lower mass bound} nor \eqref{eq:unit ball lower mass bound} and $U$ is unbounded, it may not be possible to cover $U$ with a bounded-overlap collection of balls $\{B_i\}_{i}$ for which $\mu(B_i\cap U)/\mu(B_i)$ is bounded away from one, as described in Theorem~\ref{MFBoxing}.  For example, let $X\subset\R^2$ be the completion of the region of the upper half-plane bounded between the $x$-axis and the graph of $y=1/x^2$, $x>0$, and let $U=\{(x,y)\in X:x>1\}$. Equipping $X$ with the Euclidean distance and $\mu=\Leb^2|_X$, we have that $U$ is unbounded, $\mu(U)<\infty$, and neither \eqref{eq:lower mass bound} nor \eqref{eq:unit ball lower mass bound} hold.  If $(x,y)\in U$ is contained in a ball $B$ such that $\mu(B\cap U)/\mu(B)<1/2$,  and if $x$ is sufficiently large, then $B$ must contain the set $\{(1,y):0\le y\le 1\}$.  Hence, it is not possible to obtain the cover described in Theorem~\ref{MFBoxing}. 
\end{remark}

Using the $(\theta,1,1)_{BBM}$-relative isoperimetric inequality obtained in Corollary~\ref{thm: MFRiso}, we now establish the following local version of the fractional boxing inequality: 

\begin{proposition}\label{thm:loc-box-ineq}
   Let $(X,d,\mu)$ be a complete metric measure space, with $\mu$ a doubling measure supporting a $(1,1)$-Poincar\'e inequality, and let $\kappa\in (0,1)$.  Then there exist constants $\tau\geq 1$ and $ C\ge 1$ such that for any $\mu$-measurable set $U\subset X$ and $B_0=B(x_0,r_0)\subset X$  satisfying
   \[
   \frac{\mu(U\cap B_0)}{\mu(B_0)}\leq \kappa,
   \] there exists a countable collection of pairwise disjoint balls $\{B_i\}_{i\in I\subset\N}$, with $\rad( B_i)=2^{-N_i}r_0$ for some $N_i\in\N\cup\{0\}$, and a set $\mathcal{N}\subset U$ with $\mu(\mathcal{N})=0$ such that 
    \begin{align}\label{eq:local boxing ineq balls}
        U\cap B_0\subset \mathcal{N} \cup \bigcup_{i\in I}5\tau B_i \quad \text{and}\quad \bigcup_{i\in I}B_i\subset\tau B_0,
    \end{align}
    and for all $0<\theta<1$, we have 
    \begin{align}\label{eq:loc box ineq}
        \sum_{i\in I}\frac{\mu(B_i)}{r_i^{\theta}}\le C(1-\theta)P_\theta(U,\tau B_0).
    \end{align}
   Here, $\tau$ and $C$ depend only on $\kappa$, and the doubling and $(1,1)$-Poincar\'e inequality constants. 
\end{proposition}

\begin{proof}
    We first prove the proposition under the additional assumption that $(X,d,\mu)$ is a geodesic metric measure space. Applying Lemma \ref{lem:the decomposition of the ball} with $E=U$, we obtain countably many pairwise disjoint balls $B_i \subset B_0$ , $i \in \mathbb{N}$ such that 
	\begin{align*}
    &(i)\quad\rad(B_i)=2^{-N_i}r_0\text{ for some }N_i\in\N\cup\{0\},\\
	&(ii)\quad B_0\cap U\subset \mathcal{N} \cup\bigcup_{i\in\N}5B_i \\	
	&(iii)\quad \kappa/C_{\mu}^2 \le\frac{\mu(B_i\cap U)}{\mu(B_i)}\leq\kappa,
	\end{align*}
	where $\mathcal{N} \subset U$ is a set of $\mu$-measure zero. Applying Corollary \ref{thm: MFRiso} on $B_i$ with $u=\chi_U$, it follows that 	
    \begin{align*}
	\frac{\mu(B_i)}{r_i^\theta}\leq C_1(1-\theta)&\max\left\{\frac{1}{(1-\kappa/C_\mu^2)\kappa/C_\mu^2},\frac{1}{(1-\kappa)\kappa}\right\}\\
    &\times\int_{ B_i\cap U}\int_{ B_i\setminus U}\frac{1}{d(x,y)^\theta\mu(B(x,d(x,y)))}d\mu(y)d\mu(x).
	\end{align*}
	
    Since the collection $\{ B_i\}_{i\in\N}$ has uniformly bounded radii, by the $5$-covering lemma, there exists a pairwise disjoint subcollection $\{ B_i\}_{i\in I\subset\N}$ such that 
    \[
    B_0\cap U\subset\mathcal{N}\cup\bigcup_{i\in I}5 B_i.
    \]
    Therefore, by the previous inequality and pairwise disjointness, we have 
	\begin{align*}
		\sum_{i\in I}\frac{\mu(B_i)}{r_i^\theta}&\leq C(1-\theta) \sum_{i\in I}\int_{ B_i\cap U}\int_{ B_i\setminus U} \frac{1}{d(x,y)^\theta\mu(B(x,d(x,y)))}d\mu(y)d\mu(x)\\
		&\leq C(1-\theta)\int_{ B_0\cap U}\int_{ B_0\setminus U} \frac{1}{d(x,y)^\theta\mu(B(x,d(x,y)))}d\mu(y)d\mu(x),
	\end{align*}
	where $C=C_1\max\left\{\frac{1}{(1-\kappa/C_\mu^2)\kappa/C_\mu^2},\frac{1}{(1-\kappa)\kappa}\right\}.$ Thus by taking $\tau=1$, we prove this proposition when $(X,d,\mu)$ is additionally assumed to be geodesic.

    We now prove the result without the geodesic assumption.  By the discussion at the beginning of Section~\ref{FPI}, $(X,d)$ is $L$-bilipschitz equivalent to the geodesic space $(X,\til d)$. That is, for all $x,y\in X$,
    \begin{align}\label{eq:L-equi}
        L^{-1}d(x,y)\leq \til d(x,y)\leq Ld(x,y).
    \end{align}
    From \eqref{eq:L-equi}, it directly follows that for all $x\in X$ and $r>0$,
    \begin{align}\label{eq:L-ball}
        \til B(x,r)\subset B(x,Lr)\subset \til B(x,L^2r),
    \end{align}
    where $\til B(x,r)$ denotes the ball $\til B(x,r)=\{y\in X: \til d(x,y)< r\}$. For any $\mu$-measurable set $U\subset X$ and $B_0=B(x_0,r_0)\subset X$ satisfying
   \[
   \frac{\mu(U\cap B_0)}{\mu(B_0)}\leq \kappa
   \] 
     for $0<\kappa<1$, we have by \eqref{eq:L-ball} that
   \begin{align*}
   \frac{\mu(U\cap \til B(x_0,Lr_0))}{\mu(\til B(x_0,Lr_0))}&\leq\frac{\mu(U\cap B_0)+\mu(\til B(x_0,Lr_0)\setminus B_0)}{\mu(B_0)+\mu(\til B(x_0,Lr_0)\setminus B_0)}\\
   &=\frac{\mu(U\cap B_0)/\mu(B_0)+\mu(\til B(x_0,Lr_0)\setminus B_0)/\mu(B_0)}{1+\mu(\til B(x_0,Lr_0)\setminus B_0)/\mu(B_0)}\\
   &\leq\frac{\kappa+\mu(\til B(x_0,Lr_0)\setminus B_0)/\mu(B_0)}{1+\mu(\til B(x_0,Lr_0)\setminus B_0)/\mu(B_0)}\leq \frac{\kappa+CL^{2Q_d}}{1+CL^{2Q_d}}<1,
   \end{align*}
   where the last inequality follows from the fact that $f(x)=\frac{\kappa+x}{1+x}$ is an increasing function and 
   \[
   \frac{\mu(\til B(x_0,Lr_0)\setminus B_0)}{\mu(B_0)}\le\frac{\mu(L^2B_0)}{\mu(B_0)}\le CL^{2Q_d}
   \]
   by \eqref{eq:rel lower mass bound exponent}.  By the first part of the proof under the geodesic assumption, it follows that there exists a countable collection of pairwise disjoint balls $\{\til B_i:=\til B(x_i,2^{-N_i}Lr_0)\}_{i\in I\subset\N}$ and a set $\mathcal N\subset U$ with $\mu(\mathcal N)=0$ such that 
    \begin{align*}
        U\cap \til B(x_0,Lr_0)\subset\mathcal{N}\cup\bigcup_{i\in I}5\til B_i,\qquad\bigcup_{i\in I}\til B_i\subset \til B(x_0,Lr_0),
    \end{align*}
   and such that 
   \[
   \sum_{i\in I}\frac{\mu(\wtil B_i)}{\rad(\til B_i)^\theta}\le C(1-\theta)P_\theta(U,\til B(x_0,Lr_0)).
   \]
   For each $i\in I$, let $B_i:=B(x_i,2^{-N_i}r_0)$. Then by \eqref{eq:L-ball}, it follows that      
   $\{B_i\}_{i\in I}$ is pairwise disjoint, and that
   \begin{align*}
  &(i)\quad\bigcup_{i\in I}B_i\subset L^2B_0,\\
  &(ii)\quad U\cap B_0\subset\mathcal N\cup\bigcup_{i\in I}5\til B_i\subset\mathcal N\cup\bigcup_{i\in I}5 L^2B_i, and\\
  &(iii)\quad \sum_{i\in I}\frac{\mu(B_i)}{\rad(B_i)^\theta}\lesssim\sum_{i\in I}\frac{\mu(\til B_i)}{\rad(\til B_i)^\theta}\lesssim (1-\theta)P_\theta(U,\til B(x_0,Lr_0))\le(1-\theta)P_\theta(U,L^2B_0).
  \end{align*}
   Taking $\tau=L^2$ completes the proof.\qedhere
\end{proof}

We see that the local fractional boxing inequality given in Proposition~\ref{thm:loc-box-ineq} actually implies the $(\theta,q,1)_{BBM}$-Poincar\'e inequality.

\begin{corollary}\label{cor:local box implies frac PI}
   Let $(X,d,\mu)$ be a doubling metric measure space, with relative lower mass bound exponent $Q:=Q_d>1$ \eqref{eq:rel lower mass bound exponent}, and let $0<\theta<1$. Suppose that there exist constants $C\ge 1$ and $\tau\ge 1$ such that for any $\mu$-measurable set $U\subset X$ and $B_0=B(x_0,r_0)\subset X$ satisfying
   \[
   \frac{\mu(U\cap B_0)}{\mu(B_0)}\leq 1/2,
   \] there exists a countable collection of pairwise disjoint balls $\{B_i\}_{i\in I\subset\N}$, with $\rad( B_i)=2^{-N_i}r_0$ for some $N_i\in\N\cup\{0\}$, and a set $\mathcal{N}\subset U$ with $\mu(\mathcal{N})=0$ such that 
    \begin{align*}
        U\cap B_0\subset \mathcal{N} \cup \bigcup_{i\in I}5\tau B_i \quad \text{and}\quad \bigcup_{i\in I}B_i\subset \tau B_0,
    \end{align*}
   and such that 
    \begin{align*}
        \sum_{i\in I}\frac{\mu(B_i)}{r_i^{\theta}}\le C(1-\theta)P_\theta(U,\tau B_0).
    \end{align*}
    Then $(X,d,\mu)$ supports a $(\theta,q,1)_{BBM}$-Poincar\'e inequality for all $1\le q\le Q/(Q-\theta)$, with constants depending only on $C$, $\tau$, and the doubling constant.
\end{corollary}

\begin{proof}
As in the proof of Theorem \ref{thm: MFPinq}, with $q=\frac{Q}{Q-\theta}$, it suffices to estimate the sum 
    \[
    I+II:=\int_{u_{B_0}}^{\max\{m_u,u_{B_0}\}}\mu(B_0\cap\{u>\lambda\})^{1/q}\,d\lambda+\int_{\max\{m_u,u_{B_0}\}}^\infty\mu(B_0\cap\{u>\lambda\})^{1/q}\,d\lambda,
    \]
	where we may assume $u_{B_0}<m_u$.
    Similar to \eqref{eq:I-first}, we have
    \begin{align*}
        I\le 2\int_{-\infty}^{m_u}\mu(B_0\cap\{u<\lambda\})^{1/q}\,d\lambda.
    \end{align*}
    Since $\mu(B_0\cap\{u<\lambda\})\leq \mu(B_0)/2$ for $\lambda<m_u$, we apply the hypothesis of the corollary with $U=\{u<\lambda\}$ to obtain a collection of balls $\{B_i\}_{i\in I\subset\N}$ satisfying \eqref{eq:local boxing ineq balls} and \eqref{eq:loc box ineq}.   Using also \eqref{eq:rel lower mass bound exponent}, it therefore follows that
   \begin{align*}
    \mu(B_0\cap\{u<\lambda\})^{1/q}\leq \sum_{i\in I}\mu(5\tau B_i)^{1/q}&\lesssim\sum_{i\in I}\frac{\mu(B_i)}{r_i^\theta}\frac{r_i^{\theta}}{\mu(B_i)^{(q-1)/q}}\\
    &=\frac{r_0^\theta}{\mu(B_0)^{(q-1)/q}}\sum_{i\in I}\frac{\mu(B_i)}{r_i^\theta}2^{-N_i\theta}\left(\frac{\mu(B_0)}{\mu(B_i)}\right)^{(q-1)/q}\\
    &
    \lesssim\frac{r_0^\theta}{\mu(B_0)^{(q-1)/q}}\sum_{i\in I}\frac{\mu(B_i)}{r_i^\theta}\\
    &\leq
    C(1-\theta)\frac{r_0^\theta}{\mu(B_0)^{(q-1)/q}}P_\theta(\{u<\lambda\},\tau B_0).
    \end{align*}
   By the fractional coarea formula, Lemma~\ref{lem:frac coarea}, it then follows that
    \begin{align*}
        I\le C(1-\theta)\frac{r_0^\theta}{\mu(B_0)^{(q-1)/q}}\int_{\tau B_0}\int_{\tau B_0}\frac{|u(x)-u(y)|}{d(x,y)^\theta\mu(B(x,d(x,y)))}d\mu(y)d\mu(x).
    \end{align*}
    The estimation of $II$ is the same as $I$. Therefore, the proof is finished.
\end{proof}

\section{Improved and global inequalities, and lower Ahlfors regularity}\label{Fiso}

In this section, we show that the additional assumption of lower Ahlfors $Q$-regularity \eqref{eq:lower Q} allows us to establish the improved $(\theta,\frac{Q}{Q-\theta},1)_{BBM}$-Poincar\'e and relative isoperimetric inequalities, see Corollary~\ref{cor:improved inequalities} below, as well as the $(\theta,\frac{Q}{Q-\theta},1)_{BBM}$-isoperimetric and Sobolev inequalities in Corollary~\ref{cor:Global inequalities}.  These follow essentially as applications of the local and global fractional boxing inequalities given by Proposition~\ref{thm:loc-box-ineq} and Theorem~\ref{MFBoxing}, respectively.  Finally, we prove Theorem~\ref{lowerQEquiv}, showing that these inequalities are all equivalent with lower Ahlfors $Q$-regularity of the measure $\mu$.

\begin{corollary}\label{cor:improved inequalities}
    Suppose that the assumptions of Theorem  \ref{thm: MFPinq} hold and, in addition, $\mu$ satisfies the lower Ahlfors $Q$-regular condition \eqref{eq:lower Q} with $Q>1$. Then, 
    
    \begin{enumerate}
        \item[(i)] there exist $C\ge 1$ and $\tau\ge 1$, depending only on $Q$, the doubling, lower Ahlfors regularity, and $(1,1)$-Poincar\'e inequality constants, such that for all $0<\theta<1$, and for all balls $B_0\subset X$ and $u\in L^1_\loc(X)$, we have
    \begin{align*}
        \left(\int_{B_0}|u-u_{B_0}|^{\frac{Q}{Q-\theta}}d\mu\right) ^{\frac{Q-\theta}{Q}}\le C(1-\theta)\int_{\tau B_0}\int_{\tau B_0}\frac{|u(x)-u(y)|}{d(x,y)^\theta\mu(B(x,d(x,y)))}d\mu(y)d\mu(x).
    \end{align*}
    \item[(ii)]there exist $C\ge 1$ and $\tau\ge 1$, depending only on $Q$, the doubling, lower Ahlfors regularity, and $(1,1)$-Poincar\'e inequality constants, such that for all $0<\theta<1$, and for all balls $B_0\subset X$ and measurable sets $E\subset X$, we have 
    \[
    \min\{\mu(B_0\cap E),\,\mu(B_0\setminus E)\}^{\frac{Q-\theta}{Q}}\le C(1-\theta)P_\theta(E,\tau B_0).
    \]

    \end{enumerate}
    Furthermore, $(X,d,\mu)$ supports an improved $(\theta,\frac{Q}{Q-\theta},1)_{BBM}$-Poincar\'e inequality if and only if it supports an improved $(\theta,\frac{Q}{Q-\theta},1)_{BBM}$-relative isoperimetric inequality. 
    \end{corollary}

     \begin{proof}
        We note that $(ii)$ follows from applying $(i)$ with $u=\chi_E$, and $(i)$ follows from $(ii)$ by the same argument used in the proof of Theorem~\ref{rmk: PiRiequiv}.  Thus, it remains to prove (i).  The proof is similar to  the proof of Corollary~\ref{cor:local box implies frac PI}, with minor adjustments.  We include the argument for the sake of completeness.

      As in the proof of Theorem \ref{thm: MFPinq}, with $q=\frac{Q}{Q-\theta}$, it suffices to estimate the sum 
    \[
        I+II:=\int_{u_{B_0}}^{\max\{m_u,u_{B_0}\}}\mu(B_0\cap\{u>\lambda\})^{1/q}\,d\lambda+\int_{\max\{m_u,u_{B_0}\}}^\infty\mu(B_0\cap\{u>\lambda\})^{1/q}\,d\lambda.
    \]
    Similar to \eqref{eq:I-first}, we have
    \begin{align*}
        I\le\int_{u_{B_0}}^{m_u}\mu(B_0\cap\{u>\lambda\})^{1/q}\,d\lambda&\le(m_u-u_{B_0})\mu(B_0)^{1/q}\nonumber\\
        &\le\frac{2}{\mu(B_0)^{(q-1)/q}}\int_{u_{B_0}}^{m_u}\mu(B_0\cap\{u>\lambda\})d\lambda\nonumber\\
        &\le\frac{2}{\mu(B_0)^{(q-1)/q}}\int_{u_{B_0}}^{\infty}\mu(B_0\cap\{u>\lambda\})d\lambda\nonumber\\
        &=\frac{2}{\mu(B_0)^{(q-1)/q}}\int_{-\infty}^{u_{B_0}}\mu(B_0\cap\{u<\lambda\})d\lambda\nonumber\\
        &\le\frac{2}{\mu(B_0)^{(q-1)/q}}\int_{-\infty}^{m_u}\mu(B_0\cap\{u<\lambda\})d\lambda\nonumber\\
        &\le 2\int_{-\infty}^{m_u}\mu(B_0\cap\{u<\lambda\})^{1/q}\,d\lambda.
    \end{align*}
    Since $\mu(B_0\cap\{u<\lambda\})\leq \mu(B_0)/2$ for $\lambda<m_u$, we apply Proposition \ref{thm:loc-box-ineq} with $\kappa=1/2$ and $U=\{u<\lambda\}$, to obtain a collection of balls $\{B_i\}_{i\in I\subset\N}$ satisfying \eqref{eq:local boxing ineq balls} and \eqref{eq:loc box ineq}.  Therefore, it follows that
   \begin{align*}
    \mu(B_0\cap\{u<\lambda\})^{1/q}\leq \sum_{i\in I}\mu(5\tau B_i)^{1/q}&\lesssim\sum_{i\in I}\frac{\mu(B_i)}{r_i^\theta}\frac{r_i^\theta}{\mu(B_i)^{\frac{\theta}{Q}}}\\
    &\lesssim\sum_{i\in I}\frac{\mu(B_i)}{r_i^\theta}\leq
    C(1-\theta)P_\theta(\{u<\lambda\}, \tau B_0).
    \end{align*} 
    By Lemma~\ref{lem:frac coarea}, it follows that
    \begin{align*}
        I&\lesssim (1-\theta)\int_{-\infty}^{m_u}P_\theta(\{u<\lambda\}, \tau B_0)d\lambda\\
        &\lesssim(1-\theta)\int_{ \tau B_0}\int_{ \tau B_0}\frac{|u(x)-u(y)|}{d(x,y)^\theta\mu(B(x,d(x,y)))}d\mu(y)d\mu(x).
    \end{align*}
    The estimation of $II$ is the same as $I$. Therefore, the proof is finished.
    \end{proof}

\begin{proof}[Proof of Corollary \ref{cor:Global inequalities}]
    We first prove $(i)$.  Let $E\subset X$ be measurable with $\mu(E)<\infty$. From the boxing inequality Theorem \ref{MFBoxing} and Remark \ref{Remark-lower-boxing}, 
    there exists a countable collection of disjoint balls $\{B_i:=B(x_i,r_i)\}_{i\in\mathbb{N}}$ such that 
    \[
E\subset\bigcup_{i\in\mathbb{N}}B(x_i,5\tau r_i), \quad\quad\frac{1}{2C_\mu}\le\frac{\mu(B(x_i,r_i)\cap E)}{\mu(B(x_i,r_i))}<\frac{1}{2}
    \]
    for each $i\in\mathbb{N}$, and 
\begin{equation*}
        \sum_{i\in\mathbb{N}}\frac{\mu(B(x_i,r_i))}{r_i^\theta}\le C\theta(1-\theta)P_\theta(E,X).
    \end{equation*}
Hence,
\begin{align*}
    \mu(E)^{\frac{Q-\theta}{Q}}\leq \sum_{i\in\N}\mu(5\tau B_i)^{\frac{Q-\theta}{Q}}&\lesssim\sum_{i\in\N}\frac{\mu(B_i)}{r_i^\theta}\frac{r_i^\theta}{\mu(B_i)^{\frac{\theta}{Q}}}\\
    &\lesssim\sum_{i\in\N}\frac{\mu(B_i)}{r_i^\theta}\leq
    C\theta(1-\theta)P_\theta(E,X).
\end{align*}

    To prove $(ii)$, let $u\in L^{1}(X)$. For a.e.\ $\lambda>0$, we have $\mu(\{|u|>\lambda\})<\infty.$
        Applying Cavalieri's principle and computing as in \eqref{eq:Frac PI Cavalieri}, one has 
        \[
        \left(\int_{X}|u|^{\frac{Q}{Q-\theta}}d\mu\right) ^{\frac{Q-\theta}{Q}}\leq 2\int_{0}^{\infty}\mu(\{|u|>\lambda\})^{\frac{Q-\theta}{Q}}d\lambda.
        \]
        By (i) above and Lemma~\ref{lem:frac coarea}, we have that 
        \begin{align*}
            \left(\int_{X}|u|^{\frac{Q}{Q-\theta}}d\mu\right) ^{\frac{Q-\theta}{Q}}&\leq 2\int_{0}^{\infty}\mu(\{|u|>\lambda\})^{\frac{Q-\theta}{Q}}d\lambda\\
            &\lesssim\theta(1-\theta)\int_{0}^{\infty}P_\theta(\{|u|>\lambda\},X)d\lambda\\
            &\lesssim\theta(1-\theta)\int_{X}\int_{X}\frac{|u(x)-u(y)|}{d(x,y)^\theta\mu(B(x,d(x,y)))}d\mu(y)d\mu(x).
        \end{align*}

        The proof of $(ii)$ shows that $(ii)$ follows from $(i)$.  Applying $(ii)$ with $u=\chi_E$ also shows that $(i)$ follows from $(ii)$, and so these inequalities are equivalent. 
    \end{proof}

We now prove Theorem~\ref{lowerQEquiv}.  As with Theorem~\ref{thm: MFPinq}, the conditions of the theorem persist under a bilipschitz change in metric, and so it suffices to prove the result under the additional assumption that $(X,d)$ is geodesic, see the discussion at the beginning of Section~\ref{FPI}.  We restate and prove this theorem here, under this additional assumption.

\begin{theorem}
    Let  $(X,d,\mu)$ be a complete, geodesic metric measure space, with $\mu$ a doubling measure supporting a $(1,1)$-Poincar\'e inequality such that $\mu(X)=\infty$.  Then, the following conditions are equivalent.
    \begin{enumerate}
        \item[(i)] $\mu$ satisfies the lower Ahlfors $Q$-regularity condition \eqref{eq:lower Q} for $Q>1$. 
        \item[(ii)] There exists a constant $C\ge 1$,  such that for all $0<\theta<1$ and $u\in L^1(X)$, we have
\[
\left(\int_{X}|u|^{\frac{Q}{Q-\theta}}d\mu\right) ^{\frac{Q-\theta}{Q}}\le C\theta(1-\theta)\int_{X}\int_{X}\frac{|u(x)-u(y)|}{d(x,y)^\theta\mu(B(x,d(x,y)))}d\mu(y)d\mu(x).
\]
    \item[(iii)] There exists a constant $C\ge 1$ such that for all $0<\theta<1$, and all balls $B_0:=B(x_0,r_0)\subset X$ and $u\in L^1(X)$, we have
    \[
    \left(\int_{B_0}|u-u_{B_0}|^{\frac{Q}{Q-\theta}}d\mu\right) ^{\frac{Q-\theta}{Q}}\le C(1-\theta)\int_{B_0}\int_{B_0}\frac{|u(x)-u(y)|}{d(x,y)^\theta\mu(B(x,d(x,y)))}d\mu(y)d\mu(x).
    \]
    \end{enumerate}
For the implications $(ii)\implies (i)$ and $(iii)\implies (i)$, the Ahlfors $Q$-regularity constants are independent of $\theta$.  For the implications $(i)\implies (ii)$ and $(i)\implies (iii)$, the constants $C$ depend additionally on the doubling and $(1,1)$-Poincar\'e inequality constants.
\end{theorem}

\begin{proof}
Note that $(i)\Rightarrow (ii)$ follows from Corollary~\ref{cor:Global inequalities}  while $(i)\Rightarrow (iii)$ follows from Corollary~\ref{cor:improved inequalities}. We will complete the proof by showing $(ii)\Rightarrow (i)$ and $(iii)\Rightarrow (i)$.

$(ii)\Rightarrow (i)$. For a given ball $B(z,r)$ and $ r>t$, we define a function $u(x)$ on $X$ 
by  setting
\begin{equation*}
    u(x)=\begin{cases}
        1, &\text{if} ~x\in B(z,t),\\
        \frac{r-d(x,z)}{r-t},&\text{if} ~x\in B(z,r)\setminus B(z,t),\\
        0, &\text{if} ~x\in X\setminus B(z,r).
    \end{cases}
\end{equation*}
\noindent{\it Claim:} There exists a constant $C>0$ independent of $z,t,r$ such that
\[
    \int_{X}\int_{X}\frac{|u(x)-u(y)|}{d(x,y)^\theta\mu(B(x,d(x,y)))}\,d\mu(x)\,d\mu(y)
    \le
    C\left(\frac{1}{\theta}+\frac{1}{1-\theta}\right) \frac{\mu(B(z,r))}{(r-t)^\theta}.
\]
\noindent{\it Proof of claim:} To this end, notice $u=0$ on $X\setminus B(z,r)$. Then for the above double integral, we  write
\begin{align*}
\int_{X}\int_{X}&\frac{|u(x)-u(y)|}{d(x,y)^\theta\mu(B(x,d(x,y)))}\,d\mu(x)\,d\mu(y)\\
&\lesssim\int_{B(z,r)}\int_{X\setminus B(z,r)}\frac{|u(y)|}{d(x,y)^\theta\mu(B(x,d(x,y)))}\,d\mu(x)\,d\mu(y)\\
&\quad +\int_{B(z,r)}\int_{B(z,r)}\frac{|u(x)-u(y)|}{d(x,y)^\theta\mu(B(x,d(x,y)))}\,d\mu(x)\,d\mu(y)\\
&=:H_1+H_2.
\end{align*}
For $y\in B(z,r)$, we have
\begin{align*}
	\int_{X\setminus B(z,r)}&\frac{1}{d(x,y)^{\theta}\mu(B(x,d(x,y)))}\,d\mu(x)\\
	&\le \int_{X\setminus B(y,r-d(z,y))}\frac{1}{d(x,y)^{\theta}\mu(B(x,d(x,y)))}\,d\mu(x)\\
    &=\sum_{m=0}^\infty\int_{B(y,2^{m+1}(r-d(z,y)))\setminus B(y,2^m(r-d(z,y)))}\frac{1}{d(x,y)^{\theta}\mu(B(x,d(x,y)))}\,d\mu(x)\\
	&\lesssim(r-d(z,y))^{-\theta}\sum_{m=0}^{\infty}(2^{m})^{-\theta}\\
	& =(r-d(z,y))^{-\theta}\frac{1}{1-2^{-\theta}}\\
	& \lesssim \frac{1}{\theta}(r-d(z,y))^{-\theta}.
\end{align*}
Therefore
\begin{align}\label{estimate-H1}
    H_1&\lesssim \frac{1}{\theta}\int_{B(z,r)\setminus B(z,t)}\left(\frac{r-d(y,z)}{r-t}\right) \frac{1}{(r-d(z,y))^{\theta}} d\mu(y)\notag\\
    &\quad +\frac{1}{\theta}\int_{B(z,t)}\frac{1}{(r-d(z,y))^{\theta}}d\mu(y)\notag\\
    &\lesssim \frac{1}{\theta} \frac{\mu(B(z,r))}{(r-t)^\theta}.
\end{align}
To estimate $H_2$, we write
\begin{align*}
    H_2&\le\int_{B(z,r)}\int_{B(x,r-t)}\frac{|u(x)-u(y)|}{d(x,y)^{\theta}\mu(B(x,d(x,y)))}\,d\mu(y)\,d\mu(x)\\
    &\quad +\int_{B(z,r)}\int_{B(z,r)\setminus B(x,r-t)}\frac{|u(x)-u(y)|}{d(x,y)^{\theta}\mu(B(x,d(x,y)))}\,d\mu(y)\,d\mu(x)\\
    &=:H_{2,1}+H_{2,2}.
\end{align*}
It is easy to know  that $u$ is a Lipschitz function with Lipschitz constant $1/(r-t)$. Then we have 
\[
|u(x)-u(y)|\leq \frac{1}{r-t}d(x,y).
\]
Then
\begin{align}\label{estimate: inte-ball}
   \int_{B(x,r-t)}&\frac{|u(x)-u(y)|}{d(x,y)^{\theta}\mu(B(x,d(x,y)))}\,d\mu(y)\notag\\
   &\leq \frac{1}{r-t}\int_{B(x,r-t)}\frac{1}{d(x,y)^{\theta-1}\mu(B(x,d(x,y)))}\,d\mu(y)\notag\\
   &\leq \frac{1}{r-t}\int_{B(x,r-t)}\sum_{j\geq -\log_2 (r-t)}^{\infty}C_\mu 2^{-j(1-\theta)} \frac{\chi_{B(x,2^{-j+1})\setminus B(x,2^{-j})}(y)}{\mu(B(x,2^{-j+1}))}\,d\mu(y)\qquad \text{by \eqref{eq:rho hat majorize}}\notag\\
   &\lesssim \frac{(r-t)^{1-\theta}}{(r-t)(1-\theta)}.
\end{align}
Thus
\begin{align}\label{estimate-H21}
    H_{2,1}\lesssim \frac{1}{1-\theta}\frac{\mu(B(z,r))}{(r-t)^\theta}.
\end{align}
Observing $0\leq u\leq 1$ and 
\[
\int_{X\setminus B(x,r-t)}\frac{1}{d(x,y)^{\theta}\mu(B(x,d(x,y)))}\,d\mu(y)
\lesssim 
\frac{1}{\theta}(r-t)^{-\theta},
\]
we also have
\begin{align}\label{estimate-H22}
    H_{2,2}\lesssim \frac{1}{\theta}\frac{\mu(B(z,r))}{(r-t)^\theta}.
\end{align}
Combining \eqref{estimate-H1}, \eqref{estimate-H21} and \eqref{estimate-H22}, we prove this claim.

Now let $z\in X$ and $r\in (0,\infty).$ By geodesic property, there exists a $b$ such that
\[
\mu(B(z,br))=\frac{1}{2}\mu(B(z,r)).
\]
Taking $t=br$ and by the claim, we have
\begin{align*}
    \mu(B(z,br))^{\frac{Q-\theta}{Q}}&\leq \left\lVert u\right\rVert _{L^{Q/(Q-\theta)}(X)}\\
    &\lesssim \theta(1-\theta)\int_{X}\int_{X}\frac{|u(x)-u(y)|}{d(x,y)^\theta\mu(B(x,d(x,y)))}\,d\mu(x)\,d\mu(y)\\
    &\lesssim \frac{\mu(B(z,r))}{(r-br)^\theta}.
\end{align*}
We can rewrite the above inequality as
\[
r-br\lesssim \mu(B(z,r))^{\frac{1}{Q}}.
\]
We then let $b_0=1$ and choose $b_j\in (0,1)$ satisfying the following equality
\[
\mu(B(z,b_jr))=\frac{1}{2}\mu(B(z,b_{j-1}r))=\frac{1}{2^j}\mu(B(z,r)).
\]
Obviously, $b_j\rightarrow 0$ as $j\rightarrow \infty$ and therefore
\begin{align*}
    r=&\sum_{j=1}^{\infty}(b_{j-1}r-b_j r)\\
    &\lesssim \sum_{j=1}^{\infty} \mu(B(z,b_{j-1}r))^{1/Q}\\
    &\lesssim \sum_{j=1}^{\infty} 2^{-j/Q}\mu(B(z,r))^{1/Q}\\
    &\lesssim \mu(B(z,r))^{1/Q}.
\end{align*}
Now we have completed the proof of $(ii)\Rightarrow (i)$.

$(iii)\Rightarrow (i)$: the proof of $(iii)\Rightarrow (i)$ is similar to $(ii)\Rightarrow (i)$. For completeness, we include the proof. Now taking any $z\in X$ and $r\in (0,\infty).$ Letting $b_1=1$ and by geodesic property, there exists $b_0$ and $\{b_j\}_{j=2}^{\infty}$ only depending on $X$ such that
\[
\mu(B(z,b_0r))=2\mu(B(z,r))~\text{and}~\mu(B(z,b_jr))=\frac{1}{2}\mu(B(z,b_{j-1}r))=\frac{1}{2^{j-1}}\mu(B(z,r)).
\]
Let $K\in\N$ be the smallest positive integer such that
\begin{align}\label{def-CK}
    C_K:=\left(\frac{1}{2C_s}\right)^{1/s}-\left(\frac{C}{2^K}\right)^{1/Q_d}>0.
\end{align}
Here $C\ge 1$ and $Q_d>1$ are the constants from \eqref{eq:rel lower mass bound exponent}, with $Q_d$ possibly different from $Q$ as in the hypothesis of $(iii)$, and $s$ and $C_s$ are as in \eqref{eq:intro reverse doubling}.
Then, for every $j\in\N$, we have
\[
1
>
\frac{b_j-b_{j+K-1}}{b_{j-1}}
\ge
\left(\frac{1}{2C_s}\right)^{1/s}-\left(\frac{C}{2^K}\right)^{1/Q_d}>0.
\]
Fix $j\in\N$, we define a function $u(x)$ on $B(z,b_{j-1} r)$ 
by  setting
\begin{equation*}
    u(x)=\begin{cases}
        1, &\text{if} ~x\in B(z,b_{j+K-1}r),\\
        \frac{b_{j}r-d(x,z)}{(b_{j}-b_{j+K-1})r},&\text{if} ~x\in B(z,b_{j}r)\setminus B(z,b_{j+K-1}r),\\
        0, &\text{if} ~x\in B(z,b_{j-1}r)\setminus B(z,b_{j}r).
    \end{cases}
\end{equation*}
\noindent{\it Claim:} There exists a constant $C>0$ independent of $z, \{b_j\}, r$ and $\theta$ such that
\[
    \int_{B(z,b_{j-1} r)}\int_{B(z,b_{j-1} r)}\frac{|u(x)-u(y)|}{d(x,y)^\theta\mu(B(x,d(x,y)))}\,d\mu(y)\,d\mu(x)
    \le
    C\left(\frac{1}{1-\theta}\right) \frac{\mu(B(z,b_jr))}{(b_{j}r-b_{j+K-1}r)^\theta}.
\]

\noindent{\it Proof of claim:}
To this end, we also rewrite
\begin{align*}
    \int_{B(z,b_{j-1} r)}&\int_{B(z,b_{j-1} r)}\frac{|u(x)-u(y)|}{d(x,y)^\theta\mu(B(x,d(x,y)))}\,d\mu(y)\,d\mu(x)\\
    &=\int_{B(z,b_{j-1} r)}\int_{B(z,b_{j-1}r)\cap B(x, (b_{j}-b_{j+K-1})r)}\frac{|u(x)-u(y)|}{d(x,y)^\theta\mu(B(x,d(x,y)))}\,d\mu(y)\,d\mu(x)\\
    &\quad +\int_{B(z,b_{j-1} r)}\int_{B(z,b_{j-1} r)\setminus B(x, (b_{j}-b_{j+K-1})r)}\frac{|u(x)-u(y)|}{d(x,y)^{\theta}\mu(B(x,d(x,y)))}\,d\mu(y)\,d\mu(x)\\
    &=:L_{2,1}+L_{2,2}.
\end{align*}
Obviously, $u$ is a Lipschitz function with Lipschitz constant $1/((b_{j}-b_{j+K-1})r)$. Then, by \eqref{estimate: inte-ball}, we also have
\begin{align*}
   &\int_{B(z,b_{j-1}r)\cap B(x,(b_{j}-b_{j+K-1})r)}\frac{|u(x)-u(y)|}{d(x,y)^{\theta}\mu(B(x,d(x,y)))}\,d\mu(y)\\
   &\qquad\leq \frac{1}{(b_{j}-b_{j+K-1})r}{\int_{B(z,b_{j-1}r)\cap B(x,(b_{j}-b_{j+K-1})r)}\frac{1}{d(x,y)^{\theta-1}\mu(B(x,d(x,y)))}\,d\mu(y)}\\
   &\qquad\leq \frac{1}{(b_{j}-b_{j+K-1})r}{\int_{B(x,(b_{j}-b_{j+K-1})r)}\frac{1}{d(x,y)^{\theta-1}\mu(B(x,d(x,y)))}\,d\mu(y)}\\
   &\qquad\lesssim \frac{((b_{j}-b_{j+K-1})r)^{1-\theta}}{((b_{j}-b_{j+K-1})r)(1-\theta)}.
\end{align*}
Thus
\begin{align}\label{estimate-L21}
    L_{2,1}\lesssim \frac{1}{1-\theta}\frac{\mu(B(z,b_{j-1}r))}{((b_{j}-b_{j+K-1})r)^\theta}.
\end{align}
Since  $u$ is a Lipschitz function with Lipschitz constant $1/((b_{j}-b_{j+K-1})r)$, similar to
\eqref{estimate: inte-ball} and recalling \eqref{def-CK}, we can obtain

\begin{align*}
\frac{1}{(b_{j}-b_{j+K-1})r}&\int_{B(x,2b_{j-1}r)\setminus B(x,(b_{j}-b_{j+K-1})r)}\frac{1}{d(x,y)^{\theta-1}\mu(B(x,d(x,y)))}\,d\mu(y)\\
&\lesssim 
\frac{1}{(1-\theta)((b_{j}-b_{j+K-1})r)}\left((2b_{j-1}r)^{1-\theta}-((b_{j}-b_{j+K-1})r)^{1-\theta}\right)\\
&\le\frac{1}{(1-\theta)((b_{j}-b_{j+K-1})r)}((2C_K^{-1})^{1-\theta}-1)((b_{j}-b_{j+K-1})r)^{1-\theta}.
\end{align*}
Therefore, we also have
\begin{align}\label{estimate-L22}
    L_{2,2}\lesssim \frac{1}{1-\theta}\frac{\mu(B(z,b_{j-1}r))}{((b_{j}-b_{j+K-1})r)^\theta}.
\end{align}
Combining \eqref{estimate-L21} and \eqref{estimate-L22}, the proof of the claim is finished.

Notice that if $u_{B(z,b_{j-1}r)}\ge 1/2$, then  $|u-u_{B(z,b_{j-1}r)}|\geq 1/2$ on $B(z,b_{j-1}r)\setminus B(z,b_{j}r)$, and if $u_{B(z,b_{j-1}r)}<1/2$, then $|u-u_{B(z,b_{j-1}r)}|\geq 1/2$ on $B(z,b_{j+K-1}r).$ Thus by the claim, we can obtain
\begin{align*}
\left(\frac{1}{2^{K}}\right)^{\frac{Q-\theta}{Q}} \mu(B(z,b_j r))^{\frac{Q-\theta}{Q}}&\lesssim\left\lVert u-u_{B(z,b_{j-1}r)}\right\rVert _{L^{Q/(Q-\theta)}(B(z,b_{j-1}r))}\\
    &\lesssim (1-\theta)\int_{B(z,b_{j-1}r)}\int_{B(z,b_{j-1}r)}\frac{|u(x)-u(y)|}{d(x,y)^\theta\mu(B(x,d(x,y)))}\,d\mu(y)\,d\mu(x)\\
    &\lesssim \frac{\mu(B(z,b_jr))}{(b_{j}r-b_{j+K-1}r)^\theta}.
\end{align*}
It follows that
\[
b_jr-b_{j+K-1}r\lesssim \mu(B(z,b_jr))^{\frac{1}{Q}}.
\]
The rest of the proof is the same as in $(ii)\Rightarrow (i)$ and we completed the proof of $(iii)\Rightarrow (i)$.
\end{proof}

\noindent{\bf Declarations:} The authors have no competing interests to declare that are relevant to the content of this article. No data was used in this study.
\vskip.2cm
\noindent{\bf Declaration of AI usage:} No AI tools were used in the initial version submitted to the journal.  During the revision process, however, the authors used ChatGPT Pro (GPT-6) to identify typos, and instances where the presentation, hypotheses, or logical arguments needed further clarification. The AI tool was only used in identification of these instances; all corrections and modifications were completed by the authors.

\noindent J.K.: Department of Mathematical Sciences, P.O.~Box 210025, University of Cincinnati, Cincinnati, OH~45221-0025, U.S.A.\\
\noindent E-mail: {\tt klinejp@ucmail.uc.edu}
\vskip.2cm
\noindent P.L.:  Academy of Mathematics and Systems Science, 
        Chinese Academy of Sciences, Beijing 100190, P.R. China.\\
        E-mail: {\tt panulahti@amss.ac.cn}
\vskip.2cm
\noindent J.L.: Key Laboratory of Computing and Stochastic Mathematics (Ministry of Education), School of Mathematics and Statistics, Hunan Normal University, Changsha, Hunan 410081, P.R. China.\\
E-mail: {\tt  jiang\_li@hunnu.edu.cn; jiangli.math@qq.com}
\vskip.2cm
\noindent X.Z.: Analysis on Metric Spaces Unit, Okinawa Institute of Science and Technology Graduate University,
Okinawa 904-0495, Japan.\\
E-mail: {\tt xiaodan.zhou@oist.jp}

\end{document}